\documentclass[10pt]{iopart}
\usepackage{graphicx}
\usepackage{bbm}
\usepackage{dsfont}
\usepackage[mathscr]{euscript}

\usepackage[utf8]{inputenc} % allow utf-8 input
\usepackage[T1]{fontenc}    % use 8-bit T1 fonts

\usepackage{setstack,graphicx,color,iopams,bm,dirtytalk,float,bm, comment, algpseudocode, algorithm, url}
\newcommand{\rmr}{{\rm{r}}}
\newcommand{\prox}{{\rm prox}}

\newcommand{\bx}{\bX}

\newcommand{\bX}{\mathbf{X}}
\newcommand{\bT}{\mathbf{T}}

\begin{document}

\title[Observable asymptotics of regularized Cox regression models]{Observable asymptotics of regularized Cox regression models with standard Gaussian designs: a statistical mechanics approach} 

\author{Emanuele Massa$^{\dag, \S}$ and Anthony CC Coolen$^{\dag \star}$}
\address{$\dag$  Donders Institute, Faculty of Science, Radboud University, 6525AJ Nijmegen, The Netherlands\\
$\S$ Saddle Point Science Europe, UBC Building, Mercator Science Park,  
6525EC Nijmegen, The Netherlands\\
$\star$ Saddle Point Science, Clementhorpe, York YO23 1AN, UK  
}
\ead{emanuele.massa@donders.ru.nl, a.coolen@science.ru.nl}

\begin{abstract} 
We study the asymptotic behaviour of the Regularized Maximum Partial Likelihood Estimator (RMPLE) in the proportional limit, considering an arbitrary convex regularizer and assuming that the covariates  $\mathbf{X}_i\in\mathbb{R}^{p}$  follow  a multivariate Gaussian law with covariance $\bm{I}_p/p$ for each  $i=1, \dots, n$. In order to efficiently compute the estimator under investigation, we propose a modified Approximate Message Passing (AMP) algorithm, that we name COX-AMP,  and compare its performance with the Coordinate-wise Descent (CD) algorithm, which is taken as reference. By means of the Replica method, we derive a set of six  Replica Symmetric (RS) equations that we show to correctly describe the average behaviour of the estimators when the sample size and the number of covariates is large and commensurate. These equations cannot be solved in practice, as the data generating process (that we are trying to estimate) is not known. However, the update equations of COX-AMP suggest the construction of a local field that can in turn be used to accurately estimate all the RS order parameters of the theory \emph{solely from the data}, \emph{without} actually solving the RS equations. We emphasize that this approach can be applied when the estimator is computed via any method and is not restricted to COX-AMP. 

Once the RS order parameters are estimated, we have access to the amount of signal and noise in the RMPLE, but also its generalization error, directly from the data. 
Although we focus on the Partial Likelihood objective, we envisage broader application of the methodology proposed here, for instance to GLMs with nuisance parameters, which include some non- proportional hazards models, e.g. Accelerated Failure Time models.
\end{abstract}

\section{Introduction}
\label{sec : introduction}
The Cox semi-parameteric model is ubiquitous in survival analysis, i.e. the analysis of time to event data, which arises, for instance, in medical studies.
In many modern applications the number of covariates ($p$) is comparable, or even considerably larger than, the sample size ($n$). An example is when genetic variables are included in a survival model in medicine \cite{Bradic_2011, Bradic_2021, Buhlmann_2011}. In this setting, the estimator obtained by maximization of the negative logarithm of the Cox partial likelihood  \cite{Cox_1972}
%\begin{equation}
% \label{def : cox_lpl}
%     \mathcal{PL}_n(\bbeta) =  \sum_{i=1}^n \Delta_i\Big\{\log\Big(\frac{1}{n}\sum_{j=1}^n \Theta(T_j - T_i)\rme^{\mathbf{X}_j'\bbeta}\Big) -\mathbf{X}_i'\bbeta \Big\}  \ .
% \end{equation}
does not exist for sufficiently large values of the ratio $\zeta := p / n$, \cite{Zhang_2022, Massa_22}. In order to \say{fit} the model, a regularization is required.
% , i.e.
% \begin{equation}
% \label{def : ppl}
%     \hspace{-0.5cm}\mathcal{PL}_n(\bbeta) = \sum_{i=1}^n \Delta_i\Big\{\log\Big(\frac{1}{n} \sum_{j=1}^n \Theta(T_j - T_i)\rme^{\mathbf{X}_j'\bbeta}\Big) -\mathbf{X}_i'\bbeta \Big\} + r(\bbeta) \ .
% \end{equation}
%The asymptotics statistics of the estimator above can be studied under different scaling regimes for the number of covariated included in the model ($p$). %Traditional asymptotic theory for regression was developed under the \say{classical} regime, $p =  O(1)$ and is not applicable when $p > n$.  
%This setting with the number of covariates ($p$) comparable, or even considerably larger than the sample size ($n$) is not rare in, for instance, medical applications, especially when genetic variables are included in the model \cite{Bradic_2011, Bradic_2021, Buhlmann_2011}. 

When the goal of the analysis is to obtain a more \say{parsimonious} model, sparsity inducing regularizations like e.g. Lasso \cite{Tishbirani_1996}, Elastic Net\cite{Zhou_2005}, MCP \cite{Zhang_2010} and SCAD \cite{Fan_2001} to name a few, are widely adopted. By shrinking some of the coefficients exactly to zero, these methods effectively perform variable selection, and thus control the estimation error. 
To model this situation, statisticians proposed and studied the \say{ultra high} dimensional regime  $p = O(\exp(n))$, complemented by the hypothesis that the number of active covariates (the \say{relevant} ones), $s$, is \say{small} compared to $p$ (and $n$), i.e. $s = O(n^{\delta})$ for $\delta \in (0,1)$.   
Within this regime, the sparsity inducing regularizations mentioned above, have been thoroughly investigated in the statistical literature for Generalized Linear Models \cite{Van_de_geer_2008, Buhlmann_2011} and for the Cox model \cite{Bradic_2011, Kong_2014, Bradic_2021}.
In particular, \cite{Bradic_2011} showed that, with the Lasso and folded concave regularizations (e.g. SCAD, LASSO), the Regularized Maximum Partial Likelihood Estimator (RMPLE) is equal to a (biased) oracle estimator which knows the underlying (sparse) true model, with overwhelming probability as the sample size $n$ diverges. Furthermore, the Mean Square Error of the RMPLE converges in probability to zero with the sample size $n$, and they quantify the rates of this convergence for the different regularizations.  

In this manuscript, we study the RMPLE  under a different scaling regime, namely the \say{proportional} regime where the number of observation $n$, the number of covariates $p$ and, importantly, the number of active covariates $s$, diverge proportionally, i.e. $p = \zeta n$, $s = \nu p$. In practice, it is difficult to establish which of the two regimes is more \say{realistic} and all we can do is to quantify the uncertainty and the bias of the RMPLE under hypothetical assumptions over the scaling regime of $p$ and $s$ in $n$. These quantitatively differ in the two scaling mentioned.

We highlight in the following the main findings of our study.

\paragraph{\textbf{Sharp asymptotics in the proportional regime via the replica method}}
We study the asymptotic behaviour of the RMPLE by means of the replica method under two assumptions, namely
\begin{itemize}
    \item[A1] \textbf{convex separable regularizer} 
    \item[A2] \textbf{Gaussian covariates} : $\mathbf{X}_i\sim \mathcal{N}(\bm{0}, \bm{I}_p/p)$ for all subjects $i=1,\dots,n$ \ .
\end{itemize}
Assumption A1 is required in order to avoid Replica Symmetry Breaking (RSB), which has been observed in presence of a non-convex regularization \cite{Sakata_2023}. Assumption A2 is needed to carry out the replica calculation. However, it is conjectured that the results obtained under A2 should be valid also for different distributions of the covariates, because of universality. We discuss more on this in the conclusion section.

Under the Replica Symmetric (RS) ansatz, we derive six coupled non-linear RS equations (for six order parameters), whose solution grants access to (the expectation value of) virtually any metric of interest, e.g. the Mean Squared Error per component for the estimator $\hat{\bbeta}_n$, or the fraction of false positives and negatives in support recovery (as originally noted in \cite{Salehi_2019} for logistic regression). The predictions of the theory are checked against the results of simulations. 

We emphasize that the Replica method has already been applied in the recent past to study the behaviour of estimators for the Generalized Linear Models \cite{Kabashima_2009, Takahashi_2018, Okajima_2023, Loureiro_2021, Loureiro_2022} and  for the  Cox semi-parametric model, with no regularization \cite{Massa_24, Coolen_17} or ridge regularization \cite{Coolen_2020, Sheikh_2019}.
Our results show, in line with the previously mentioned papers, that, in the proportional regime, the Mean Squared Error per component of the RMPLE $\hat{\bbeta}_n$ converges (in probability) to a strictly positive constant as the sample size ($n$) increases, and that this constant can be computed once the RS equations are solved. So, as expected, the results of \cite{Bradic_2011, Kong_2014} are not valid in this scaling regime. 

\paragraph{\textbf{Estimation of the Replica Symmetric order parameters from the data in the proportional regime }}  
Under the assumptions A1-A2, we propose and test a generalization of the Approximate Message Passing algorithm \cite{Donoho_2009, Rangan_2010} for the Cox model, which we refer to as COX-AMP, in order to compute the RMPLE efficiently in this regime. We show by numerical simulations that the resulting algorithm leads to estimators that are very close (in L2 distance) to the estimators obtained via the Coordinate-wise Descent (CD) algorithm \cite{Zou_05}, which is taken as a reference.
The update equations of COX - AMP suggest the construction of a \say{local field}, that we  use to estimate the order parameters of the theory from the data, i.e. without the need of actually solving the RS equations. We show via numerical simulations, that these estimators, obtained without the knowledge of the data generating process (besides what is assumed above), perform extremely well, i.e. they are very close to the numerical solution of the RS equations, which conversely requires the knowledge of the data generating process. Always by numerical simulations, we show that the \say{local field} might be computed directly also from the estimator obtained by the CD algorithm (or other equivalent ones). Our approach is a generalization of what is done by means of the Expectation Consistent method \cite{Opper_2005} by the authors of \cite{Takahashi_2018}. 
However, to the best of our knowledge, their methodology cannot be directly applied here.

We highlight that estimating the RS order parameter from the data is preferable to the direct solution of the RS equations, since the knowledge of the data generating process is not required. This is fundamental in order to use in practice the knowledge extracted from the RS equations, since with real data the only available quantities are the data and the estimators (computed from the data). To make this point clear, we stress that knowledge of the RS order parameters allows computing the amount of signal and noise in the estimator, as quantified by the overlaps of the theory which, in turn, allow the computation of the generalization error and, eventually, debiasing when the signal strength is estimated with a modified version of the algorithm proposed in \cite{Massa_24}. These applications will be the subject of a forthcoming manuscript.

\section{Background on survival analysis and the Cox regression model}

A time-to-event data point is a measurement of the amount of time that is passed between a predefined time origin and the event of interest. To be concrete, the predefined time origin might be the time at which a patient is admitted in a medical study, and the event of interest can be death, or cancer relapse.
Time to event data sets are often incomplete. For some subjects we might only know that the event  of interest (e.g. death) occurred after (before) a certain time point, in this case we speak of right (left) censoring. An example for right censoring is that a patient leaves the study since she/he changes hospital,  whilst for left censoring an example is the time of the diagnosis (a sickness can be diagnosed only after being contracted). Other types of censoring exist, for instance when it is known that the event under investigation verified within a specific interval, e.g. cancer relapse. Focusing on the case of right censored data for simplicity, time to event data comprise observations reporting: i) the event indicator $\Delta_i$, equal to $1$ if the subject experienced the event and $0$ otherwise, ii) the observed event time $T_i$, and iii) the covariate vector $\mathbf{X}_i = (X_{1},\dots,X_{p})\in \mathbb{R}^p$, i.e. the list of characteristics of the subject.

A possible model for right censored time to event data is that the event indicator $\Delta_i$ and event time $T_i$ for the i-th subject are generated as 
\begin{equation}
\label{def : gen_mechanism}
    \Delta_i = \Theta(C_i-Y_i), \quad T_i = \min\{Y_i, C_i\}
\end{equation}
where $C \sim f_C$ is the (latent) random variable censoring time and $Y_i$ is the (latent) censoring time. 
In words, the subject experiences the event during the study if she/ he does not drop out of the study first.
A regression model can then be constructed by letting the conditional density of the latent event time $Y_i$ depend on the features $X_i$, i.e.
\begin{equation}
\label{def : cond_dens_y}
    Y_i|\mathbf{X}_i \sim f(.|\mathbf{X}_i) = \lambda(.|\mathbf{X}_i) \exp\big\{-\Lambda(.|\mathbf{X}_i)\big\}  \ ,
\end{equation}
where $\lambda(.|\mathbf{X}_i)$ is the hazard rate, defined as the probability per unit time for subject $i$ to experience the event at time $t$ given that she/ he did not experience the event until time $t$ and its covariates vector equals 
$\mathbf{X}_i$. Mathematically,
\begin{equation}
    \lambda(t|\mathbf{X}_i) \rmd t := P\Big[t< Y_i < t + dt \Big| Y_i> t, \mathbf{X}_i\Big] \ .
\end{equation}
Under the assumption that $Y_i$ is independent of $C_i$ given $X_i$ (usually referred to as uninformative censoring), the (unknown) generative model is 
\begin{equation}
    f(\Delta, t|\mathbf{X}_i) = f(t|\mathbf{X}_i)^{\Delta} S_C(t|\bX_i)^{\Delta} \times S(t|\bX_i)^{1-\Delta} f_c(t|\bX_i)^{1-\Delta} \ ,
\end{equation}
where $f_c\ , S_C(x) = \int_x^{\infty} f_c(x) \rmd x $ are, respectively, the density and survival functions for the censoring risk, while 
\begin{equation}
    S(t|\mathbf{X}_i) := \int_t^{\infty} f(s|\mathbf{X}_i'\bbeta) \rmd s = \exp\big\{-\Lambda(t|\mathbf{X}_i)\big\} \ ,
\end{equation}
is the survival function of the primary risk under investigation.

\paragraph{ \textbf{Proportional hazards models}}
The class of Proportional Hazard (PH) models is one of the most widely adopted in medical applications \cite{Kalbfleisch_2011, Harrell_2001, Cox_1972}. This assumes the following (semi-) parametrization
\begin{equation}
\label{assumption : ph}
    \lambda(t|\mathbf{X}_i) = \lambda(t) \exp(\bX_i\bbeta)\ .  
\end{equation}
The above implies that
\begin{eqnarray}
\label{def : surv_ph}
    S(t|\mathbf{X}_i'\bbeta) &:=& \exp\{-\Lambda(t)\exp(\mathbf{X}_i'\bbeta)\},  \\ 
\label{def : dens_ph}
    f(t|\mathbf{X}_i'\bbeta) &:=&  \lambda(t)\exp \{ \mathbf{X}_i'\bbeta - \Lambda(t)\exp(\mathbf{X}_i'\bbeta)\}, 
\end{eqnarray}
where $\Lambda(t)$ is the base cumulative hazard, and $\lambda(t) := \rmd \Lambda(t) / \rmd t$ is the base hazard rate for the primary risk. The latter is the probability per unit time of the \say{average} subject, i.e. $\bX_i = \bm{0}$ when the covariates are centered, to experience the event in the infinitesimal interval $(t, t + \rmd t)$, given she/he did not experience the event until $t$.  
Neglecting terms that do not depend on the parameters of interest, the  minus log-likelihood for a data-set of i.i.d observations $\{(\Delta_i,T_i,\mathbf{X}_i)\}_{i=1}^n$ reads
\begin{eqnarray}
\label{def : ll_ph}
    \mathcal{L}_n\big(\bbeta,\lambda\big) &=& \sum_{i=1}^n \Big\{ \Lambda(T_i)\rme^{\mathbf{X}_i'\bbeta} - \Delta_i\big(\log \lambda(T_i) +\mathbf{X}_i'\bbeta\big)\Big\}\\
    &=& \sum_{i=1}^n \Big\{g\big(\mathbf{X}_i'\bbeta, \Lambda(T_i), \Delta_i\big)- \Delta_i \log \lambda(T_i)\Big\} \ ,
\end{eqnarray}
where we introduced 
\begin{equation}
    g(x,y,z) :=  \exp(x)y - z x\ , 
\end{equation}
in order to highlight the part of the likelihood function that depends on the linear predictors $\bX_i'\bbeta$.
By assuming a parametric form for the cumulative hazard rate $\Lambda$, one obtains a parametric model, see \cite{Kalbfleisch_2011, Harrell_2001} for an introduction to parametric proportional hazards survival models. In that case, the model becomes a parametric GLM with additional nuisance parameters, i.e. the parameters of the base hazard rate. 
The Cox semi-parametric model \cite{Cox_1972, Harrell_2001} does not make any assumption on the shape of the hazard function, leaving the latter unspecified, as a consequence it is widely adopted in applications.  

Since the assumption (\ref{assumption : ph}) implies that the ratio of the cumulative hazard of two subjects is independent of time, several methods have been proposed to test the validity of the PH assumption, see \cite{Harrell_2001} for a more detailed discussion. When the latter is violated, several options are available. For instance, assuming that the base hazard rate varies among subgroups in the population, e.g. between males and females, this is called stratification of the model. In this case, the PH assumption does not hold at the population level, but does hold within the subgroups.
It is also possible to assume time varying coefficients within the PH parametrization, so that the PH assumption can be violated. There also exist other survival analysis models that make different assumptions rather than the PH one, e.g. Accelerated Failure Time (AFT) models, see \cite{Harrell_2001, Kalbfleisch_2011}.

\paragraph{\textbf{Cox partial likelihood}}
We now give a \say{quick and dirty} derivation of the Cox Partial Likelihood \cite{Cox_1972}.
By first minimizing the functional likelihood (\ref{def : ll_ph}) with respect to the function $\lambda$, see \cite{Massa_24}, one obtains the so-called Breslow estimator \cite{Breslow_72, Friedman_82}:
\begin{equation}
\label{def : breslow}
     \hat{\lambda}_{n}(t,\bbeta) =  \sum_{k=1}^n  \frac{\Delta_k\delta(t-t_k)}{\sum_{j=1}^n \Theta(t_j - t_k)\rme^{\mathbf{X}_j' \bbeta} } \ ,
\end{equation}
where $\Theta(x)$ denotes the Heaviside step-function, equal to one if $x \geq 0$ and zero else, and $\delta(x)$ the Dirac's delta distribution. Substituting in (\ref{def : ll_ph}) and disregarding terms which do not depend on $\bbeta$ it is possible to obtain the logarithm of the Cox partial likelihood 
\begin{equation}
\label{def : cox_lpl}
    \mathcal{PL}_n(\bbeta) =  \sum_{i=1}^n \Delta_i\Big\{\log\Big(\frac{1}{n}\sum_{j=1}^n \Theta(T_j - T_i)\rme^{\mathbf{X}_j'\bbeta}\Big) -\mathbf{X}_i'\bbeta \Big\}  \ .
\end{equation}
The function above is convex in $\bbeta$ and thus can be minimized using any convex optimization algorithm.
Optimizing the Cox partial likelihood, one obtains an estimator $\hat{\bbeta}_n$ and, consequently, an estimator for the cumulative hazard function of the primary risk via the Nelson-Aalen estimator
\begin{equation}
\label{def: NA_est}
    \hat{\Lambda}_n(t) :=  \sum_{i=1}^n \frac{\Delta_i \Theta(t-T_i)}{\sum_{j=1}^n \Theta(T_j-T_i)\rme^{\mathbf{X}_j'\hat{\bbeta}_n}}  
\end{equation}
which is obtained via integration of (\ref{def : breslow}) with respect to $t$.

When $p$ is sufficiently large, the Maximum Partial Likelihood estimator, obtained by maximization of (\ref{def : cox_lpl}),  does not exist \cite{Zhang_2022, Chen_2009}, and the inclusion of a regularization term $\rmr$ is required in order to define a well-posed optimization problem. With this in mind, we re-define
\begin{equation}
\label{def : ppl}
    \hspace{-0.5cm}\mathcal{PL}_n(\bbeta) = \sum_{i=1}^n \Delta_i\Big\{\log\Big(\frac{1}{n} \sum_{j=1}^n \Theta(T_j - T_i)\rme^{\mathbf{X}_j'\bbeta}\Big) -\mathbf{X}_i'\bbeta \Big\} + r(\bbeta) \ .
\end{equation}

\section{Computing the  Regularized Maximum Partial Likelihood Estimator (RMPLE)}
\label{sec : AMP_CD}

When the covariates are not correlated and Gaussian with mean zero, a variant of the Approximate Message Passing (AMP) algorithm \cite{Donoho_2009} has been proposed to compute the Maximum A Posteriori estimator for Generalized Linear Models \cite{Rangan_2010}. We propose to use a modified version of the Generalized-AMP, in order to compute the Penalized Partial Likelihood estimator. We call the resulting algorithm COX-AMP, see Algorithm \ref{alg : amp_cox}, we give a derivation in \ref{app : amp_cox}.

\begin{algorithm}
\begin{algorithmic}[1]
\Require{$\hat{\bbeta}^{0}, \ \bm{\xi}^{0}, \ \tau^0, \ \hat{\tau}^0$}
\State $\zeta \gets p / n$,  tol $\gets 1.0e-8$, max epochs $\gets 1000$
\State flag $=$ True,  err $= \infty$, t $= 0$ 
\While{err $\geq$ tol and flag}
\State $t \gets t + 1$
\State $\hat{\Lambda}_n^{t}(\bT) \gets {\rm NA}(\bT, \prox_{g(., \hat{\Lambda}^{t-1}(\bT),\bm{\Delta})}(\bxi^{t-1}, \tau^{t-1}))$
\State  err $\gets \ \max(|\Lambda(\bT)^{t+1} - \Lambda(\bT)^{t-1}|)^2 $
\State $\bm{\xi}^{t} \gets \mathbf{X}\hat{\bbeta}^{t-1} + \tau^{t-1} \dot{\mathcal{M}}_{g(., \hat{\Lambda}^t(\bT),\bm{\Delta})}(\bm{\xi}^{t-1}, \tau^{t-1})$
\State err $\gets $ err + $\max(|\bxi^{t} -\bxi^{t-1}|)^2$
\State $\hat{\tau}^{t} \gets \zeta / \Big\langle \ddot{\mathcal{M}}_{g(., \hat{\Lambda}^t(\bT),\bm{\Delta})}(\bm{\xi}^t, \tau^{t-1}) \Big\rangle $ 
\State err $\gets $ err + $(\hat{\tau}^{t} -\hat{\tau}^{t-1})^2$
\State $\bm{\psi}^t \gets  \hat{\bbeta}^{t-1} - \hat{\tau}^t \mathbf{X}'\dot{\mathcal{M}}_{g(., \hat{\Lambda}^t(\bT),\bm{\Delta})}(\bm{\xi}^t, \tau^{t-1})$
\State $\hat{\bbeta}^{t}\gets \prox_{r(.)}(\bm{\psi}^t, \hat{\tau}^t)$
\State err $\gets $ err + $\max(|\hat{\bbeta}^{t} -\hat{\bbeta}^{t-1}|)^2 $
\State $\tau^{t} \gets    \hat{\tau}^t  \ \Big\langle \prox'_{r(.)}(\bm{\psi}^{t},\hat{\tau}^t) \Big\rangle$
\State err $\gets $ err + $(\tau^{t} -\tau^{t-1})^2$
\State err $\gets \sqrt{{\rm err}}$
\If{t $\geq$ max epochs}
\State{flag $=$ False}
\EndIf
\EndWhile
\end{algorithmic}
\caption{COX-AMP algorithm for MAP in Cox regression}
\label{alg : amp_cox}
\end{algorithm}

A strategy that is widely applied in practical applications is, for instance, Coordinate-wise Descent (CD)\cite{Simon_11}. This algorithm is easy to implement and has the advantage of being applicable also when the regularizer is not differentiable, e.g. for the Lasso. Furthermore, it is derived without explicitly using the distribution of the covariates, and is, thus, more generally applicable than COX-AMP, which conversely relies on the assumption of uncorrelated Gaussian covariates with variance of order $1/p$. For this reason, we take the CD algorithm as standard reference to benchmark the behaviour of COX-AMP.

We denote
\begin{equation}
     {\rm NA}(T_i, \mathbf{X}\hat{\bbeta}_n) := \sum_{l=1}^n \frac{\Delta_i \Theta(T_i-T_l)}{\sum_{j=1}^n \Theta(T_j-T_l)\rme^{\mathbf{X}_j'\hat{\bbeta}_n}} \ , 
\end{equation}
and 
\begin{equation}
    {\rm NA}(\bT, \bX\hat{\bbeta}_n) = \Big({\rm NA}(T_1, \mathbf{X}_1'\hat{\bbeta}_n), \dots, {\rm NA}(T_n, \mathbf{X}_n'\hat{\bbeta}_n)\Big) \ .
\end{equation}
A variant of the Coordinate-wise Descent algorithm of \cite{Zou_05, Simon_11} can be summarized as in Algorithm \ref{alg : cd_cox}.
\begin{algorithm}
\begin{algorithmic}[1]
\Require $\hat{\bbeta}^{0}$
\State $\zeta \gets p / n$,  tol $\gets 1.0e-8$, max epochs $\gets 100$
\State flag $\gets$ True,  err $\gets \infty$, t $\gets 0$ 
\State $\bm{\Lambda}^{0} \gets {\rm NA}(\bT, \bX\hat{\bbeta}^{0})$
\While{err $\geq$ tol}
\State $\mathbf{s}(\bbeta^{t},\Lambda^{t}) \gets \Big\{\Lambda^{t}(\bT)\rme^{\mathbf{X}\bbeta^{t}} - \bm{\Delta}\Big\}\mathbf{X}$
\State $\mathbf{M}(\bbeta^{t},\Lambda^{t}) \gets \bX '\Big\{\Lambda^t(\bT)\rme^{\mathbf{X}\bbeta^t} \Big\} \mathbf{X} $
\State $\bvarphi^{t} \gets \bbeta^{t}$
\For{$1\leq k\leq p$}
\State $\psi_k \gets \big\{\mathbf{e}_k'\mathbf{M}(\bbeta^{t}, \Lambda^{t})\big\{\bbeta^{t} -(\bm{I} - \mathbf{e}_k\mathbf{e}_k')\bvarphi^{t}\big\}- s_k(\bbeta^{t}, \Lambda^{t})\big\}/\mathbf{e}_k'\mathbf{M}(\bbeta^{t}, \Lambda^{t})\mathbf{e}_k$
\State $1/\tau_k \gets \mathbf{e}_k'\mathbf{M}(\bbeta^{t}, \Lambda^{t})\mathbf{e}_k$
\State $ \varphi_k^t \gets {\rm prox}_{r(.)}\big(\psi_k, \tau_k\big)$
\EndFor
\State $\bbeta^{t+1} \gets \bvarphi^{t}$
\State $\Lambda(\bT)^{t+1} \gets {\rm NA}(\bT, \mathbf{X}\bbeta^{t+1} ) $
\State ${\rm err} \gets \sqrt{ \max (|\bbeta^{t+1} - \bbeta^{t}|)^2 + \max(|\Lambda(\bT)^{t+1} - \Lambda(\bT)^{t}|)^2 }$
\EndWhile
\end{algorithmic}
\caption{Coordinate wise path solution for Cox model}
\label{alg : cd_cox}
\end{algorithm}
We give a brief derivation in \ref{app : cd_algorithm}.

In this manuscript, the elastic net regularization 
\begin{equation}
\label{def : elastic_net}
    \rmr(\bbeta) := \alpha  \|\bbeta\|_1 + \eta \frac{1}{2} \|\bbeta\|_2^2 \  , 
\end{equation}
will be taken as the reference convex separable regularization. It is handy to use the alternative parametrization
\begin{equation}
    \alpha = \rho \lambda, \quad \eta = \rho(1-\lambda) \ , 
\end{equation}
where $\rho$ is the strength of the regularization and $\lambda$ is the interpolation parameter between the L1 and L2 norm. That is, when $\lambda = 1.0$, the regularizer is the Lasso (L1 norm), whilst when $\lambda = 0.0$ the regularizer equals the Ridge (L2 norm). 

\begin{figure}[H]
\centering
\includegraphics[scale = 0.7]{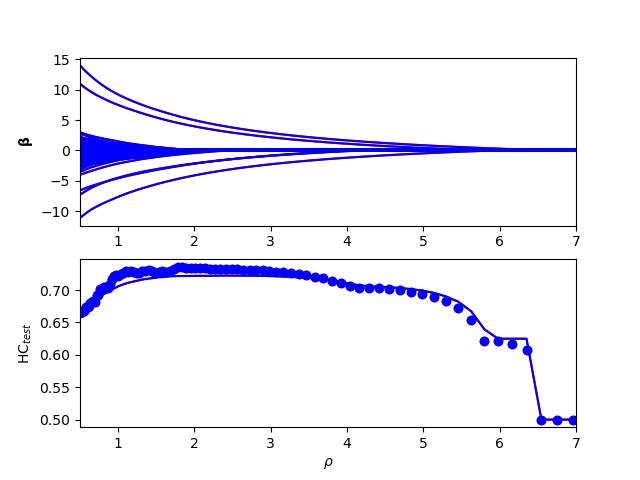}
\caption{\textbf{Elbow plot of $\hat{\bbeta}_n(\alpha)$ (above) and Harrell's C-index (below)} computed via Cox-AMP in red and Cox-CD in blue. 
The data are generated from a Log-logistic proportional hazards model $\Lambda_0(t) = \log(1+t^2/2)$ with uniform censoring between $\tau_1 = 1.0$ and $\tau_2 = 2.0$. The associations  $\bbeta_0$ are sampled as from a Gauss-Bernoulli distribution $\mathcal{P}_{\beta_0}(x) = (1-\nu) \delta(x) + \nu \exp\{-x^2/2\sigma^2\} / \sqrt{2\pi}\sigma$, with sparsity $\nu = 0.005$ and $\sigma$ adjusted to have $\|\bbeta_0\|^2_2/p = \theta^2_0 = 1.0$. The L1 ratio $\lambda = 0.75$.  The L2 relative distance between $\hat{\bbeta}_n(\alpha)$ computed via Cox-AMP and Cox-CD is of the order of $10^{-6}$ for all value of the regularization parameter displayed. The coloured dots are the estimate of the test c - index defined alter in (\ref{def : CV_C_index}), obtained via the identity (\ref{def : tilde_bxi}), obtained by the replica symmetric theory (using the training set only).}
\label{fig : elbow_plot}
\end{figure}

Numerical simulations show, see figures (\ref{fig : elbow_plot}), that the estimator computed via Cox - AMP is very close (in L2 norm) to the estimator computed via Coordinate-wise Descent. We noticed that Cox - AMP, albeit being generally faster than CD, might not always converge and, generally, requires a damped update, especially when the regularization strength is \say{small}. 
It is well known that the Maximum Partial Likelihood estimator does not exist (with probability approaching one in the asymptotic limit) past a critical threshold, $\zeta_c$, which depends on the data generating process.
In practice, the estimate $\hat{\bbeta}_n$ diverges to infinity (in L2 norm) as $\rho\rightarrow 0^+$ \cite{Chen_2009, Zhang_2022}. 
Indeed, in all our simulations, both the Cox AMP and the Cox - CD algorithms do not converge at sufficiently small regularization strength. 
However, these values of the regularization parameter are unlikely to be selected with cross validation, since, as can be seen in figure (\ref{fig : elbow_plot}), they return a suboptimal prediction error, compared to the optimal one which is obtained for larger values of the regularization strength.
We quantified the prediction accuracy by the Harrell's Concordance index,
\cite{Harrell_82}, defined for a sample of $n_{test}$ unseen subjects as follows 
\begin{equation}
    \label{test_C_index}
        \widehat{HC}_{n_{test}} =\frac{ \sum_{i=1}^{n_{test}} \Delta_i \sum_{j=1}^{n_{test}} \Theta(T_j-T_i) \Theta\Big(\mathbf{X}_j'\hat{\bbeta}_n-\mathbf{X}_i'\hat{\bbeta}_n\Big)}{ \sum_{i=1}^{n_{test}} \Delta_i \sum_{j=1}^{n_{test}} \Theta(T_j-T_i) }  \ .\nonumber
\end{equation}
This metric measures the ability of the model to separate (unseen) individuals with different risks scores. The rationale being that a good model should assign shorter survival times to subjects with higher risk score \cite{Steyenberg_2010}.
This quantity is $1$ if the model has perfect discrimination ability. Random guessing gives a c-index of $0.5$.  
The actual test C-index is availble in this case, because we simulated the data and is depicted as a continuous line (red for Cox - AMP and blue for Cox- CD, but only the blue is visible). In applications, one needs to estimate this quantity, for instance by cross validation \cite{Van_Houwelingen_1993, Le_cessie_1990}. We show in the following, see subsection \ref{subsec : application}, that the Replica Theory suggests a very fast way to estimate the test C - index, namely (\ref{def : CV_C_index}). This estimate is reported as coloured dots (again red for Cox -AMP and blue for Cox - CD, only blue is visible).  

\section{Typical behaviour of the RMPLE}
\label{section:replica_results}
In the statistical physics approach to optimization, the optimal value of the objective function (\ref{def : ppl}) is equivalent to minus the zero temperature free energy density of a fictitious physical system, i.e. 
\begin{equation}
\label{eq: laplace_identity}
    \lim_{n\rightarrow \infty} \mathbb{E}_{\mathcal{D}}\Big[\frac{1}{n}\mathcal{PL}_n(\hat{\bbeta}_n) \Big]  = \lim_{n\rightarrow \infty}  \lim_{\gamma\rightarrow \infty} \frac{1}{n \gamma }\mathbb{E}_{\mathcal{D}}\Big[ -\log \int \rme^{-\gamma \mathcal{H}_n(\bbeta|\mathcal{D})} \rmd \bbeta\Big]  , 
\end{equation}
with Hamiltonian equal to the minus penalized partial likelihood 
\begin{equation}
\label{eq : Hamiltonian}
    \mathcal{H}_n(\bbeta|\mathcal{D}) = \sum_{i=1}^n \Delta_i \Big[\log\Big(\frac{1}{n}\sum_{j=1}^n \Theta(T_j - T_i)\rme^{\mathbf{X}_i'\bbeta}\Big) - \mathbf{X}_i'\bbeta  \Big] + r(\bbeta)
\end{equation}
where we indicated with $\mathcal{D}$ the data-set, i.e. the set of tuples $\{\Delta_j,T_j,\mathbf{X}_j\}_{j=1}^n$.
The $\bbeta$s are the degrees of freedom of the system and the data-set $\mathcal{D} = \{\Delta_j,T_j,\mathbf{X}_j\}_{j=1}^n$ plays the role of the quenched disorder.

We compute the right-hand side of (\ref{eq: laplace_identity}) via the replica method under A1-A2, in the proportional regime, and assuming that the event indicator and event time are generated as in (\ref{def : gen_mechanism}), where the latent event time $Y_i$ has conditional density as in (\ref{def : cond_dens_y}).
The scaling $1/\sqrt{p}$ of the covariates is a standard assumption in virtually any previous study of high dimensional regression,  see for instance \cite{Coolen_17, Coolen_2020, Loureiro_2022, Gerace_2021} and is required to avoid introducing \say{ad hoc} entropic constants, see \cite{Massa_24}, to have well-defined $n,p\rightarrow\infty$ after the exchange with the $r \rightarrow 0$ limit, with $r$ the number of replicas. Furthermore, in practice, one might rescale the parameters of the regularizer instead. For instance, consider the elastic net regularization
\begin{equation}
    \rmr(\bbeta; \alpha, \eta) = \alpha \|\bbeta\|_1 + \frac{1}{2}\eta\|\bbeta\|_2^2 \ ,
\end{equation}
then 
\begin{eqnarray*}
    \hspace{-2cm}&&\hat{\bbeta}_n := \underset{\bbeta\in \mathbb{R}^p}{\arg\min} \Big\{\sum_{i=1}^n \Delta_i \Big[\log\Big(\frac{1}{n}\sum_{j=1}^n \Theta(T_j - T_i)\rme^{\mathbf{X}_i'\bbeta}\Big) - \mathbf{X}_i'\bbeta  \Big] + \rmr(\bbeta; \alpha, \eta)\Big\} \\
    \hspace{-2cm}&& \hspace{0.2cm}= \sqrt{p}  \ \underset{\bbeta\in \mathbb{R}^p}{\arg\min} \Big\{\sum_{i=1}^n \Delta_i \Big[\log\Big(\frac{1}{n}\sum_{j=1}^n \Theta(T_j - T_i)\rme^{'\tilde{\mathbf{X}}_i'\bbeta}\Big) - \mathbf{X}_i'\bbeta  \Big] + \rmr(\bbeta; \sqrt{p}\alpha, p \eta)\Big\} \ , 
\end{eqnarray*}
with $\tilde{\mathbf{X}}_i \sim \mathcal{N}(0,1)$.

The full replica derivation is included in \ref{app_sec : replica_computation}, and gives
\begin{equation}
\label{extr}
     \lim_{n\rightarrow \infty} \mathbb{E}_{\mathcal{D}}\Big[\frac{1}{n}\mathcal{PL}_n(\hat{\bbeta}_n) \Big] = \underset{w,v,\tau, \hat{w}, \hat{v}, \hat{\tau}}{{\rm extr}}\mathcal{F}(w,v,\tau, \hat{w}, \hat{v}, \hat{\tau}) \ .
\end{equation}
We now introduce some useful definitions that will recur throughout the manuscript and in particular to define $\mathcal{F}$. For a convex function $b:\mathbb{R}^{d} \rightarrow \mathbb{R}$, according to \cite{Rockafellar_1997}, we introduce:
\begin{itemize}
    \item  the Moureau envelope $\mathcal{M}_{b(.)}$ 
\begin{equation}
    \mathcal{M}_{b(.)}(\mathbf{x},\alpha) = \underset{\mathbf{z}\in \mathbb{R}^d}{\min}\Big\{\frac{1}{2\alpha} \|\mathbf{z} - \mathbf{x}\|^2  + b(\mathbf{z})\Big\}, 
\end{equation}
    \item the proximal mapping operator $\prox_{b(.)}$ 
\begin{equation}
\label{def : proximal}
    \prox_{b(.)}(\mathbf{x},\alpha) = \underset{\mathbf{z}\in \mathbb{R}^d}{\arg\min}\Big\{\frac{1}{2\alpha} \|\mathbf{z} - \mathbf{x}\|^2  + b(\mathbf{z})\Big\} \ . 
\end{equation}
\end{itemize}
Furthermore it is handy to define the limit of the empirical distribution of the entries of $\bbeta_0$ as 
\begin{equation}
\label{def : p_beta0}
    \mathcal{P}_{\beta_0}(x) := \lim_{p\rightarrow\infty } \frac{1}{p} \sum_{\mu=1}^p \delta(x - \mathbf{e}_{\mu}'\bbeta_0) \ ,
\end{equation}
so that a component of $\bbeta_0$ can be seen as a draw from the distribution (\ref{def : p_beta0}) above.
The zero temperature (disordered averaged) free energy in the Replica Symmetric ansatz reads 
\begin{eqnarray}
    \hspace{-1cm}\mathcal{F}(w,v,\tau,\hat{w},\hat{v},\hat{\tau})&=& -\frac{\zeta}{2\hat{\tau}}  \Big((w - \hat{w})^2 + v^2 + \hat{v}^2(1-\tau/\hat{\tau})\Big) \nonumber \\
     \hspace{-1cm}&+& \zeta \mathbb{E}_{Z, \beta_0}\Big[\mathcal{M}_{\rmr(.)}\Big(\hat{w}\frac{\beta_{0}}{\theta_0} +  \hat{v}Z, \hat{\tau}\Big)\Big]  \nonumber \\
    \hspace{-1cm}&+& \mathbb{E}_{\Delta,T,Z_0,Q}\Big[\mathcal{M}_{ g(.,\Lambda(T),\Delta)}\Big( w Z_0 + v Q , \tau\Big)\Big] + {\rm const}
\end{eqnarray}
where $Z_0,Q \sim \mathcal{N}(0,1)$, $Z_0\perp Q$, 
\begin{equation}
    \Delta, T|Z_0 \sim  f(\Delta, t|\theta_0 Z_0), \ \theta_0 := \|\bbeta_0\|/\sqrt{p}, \ Z_0,Q \sim \mathcal{N}(0,1), \ Z_0\perp Q \ , 
\end{equation}
and the function $\Lambda$ is defined at $w, v, \tau$ as 
\begin{eqnarray}
     \label{eq: self_const_Lambda_1}
     \Lambda(t) &=& \mathbb{E}_{\Delta,T,Z_0,Q}\Bigg[ \frac{ \Delta\Theta(t-T)}{\mathcal{S}(T)}\Bigg], \\
     \label{eq: self_const_Lambda_2}
     \mathcal{S}(t) &=& \mathbb{E}_{\Delta,T,Z_0,Q}\Big[\Theta(T- t)\rme^{w Z_0 + v Q + \tau \Delta - W_0\big(\tau\rme^{\tau\Delta + wZ_0 + v Q}\Lambda(T)\big)}\Big], 
\end{eqnarray}
with $W_0(.)$ the  (real branch of) Lambert W-function, i.e. the solution of $W_0(x)\exp W_0(x) = x, \ x >-1/\rme$.

The point of extremum in (\ref{extr}) satisfies the replica symmetric (RS) equations (see \ref{app_sec:rs_eqs} for explicit derivation)
\begin{eqnarray}
\label{eq : rs1}
    w  &=& \mathbb{E}_{Z, \beta_0}\Big[\beta_{0} \varphi\Big]/\theta_0\\
\label{eq : rs2}
    \hat{v}\frac{\tau}{\hat{\tau}} &=&\mathbb{E}_{Z, \beta_0}\Big[Z \varphi\Big] \\
\label{eq : rs3}
    (w^2+v^2) &=&  \mathbb{E}_{Z, \beta_0}\Big[\varphi^2\Big]\\
\label{eq : rs4} 
     \hat{w} &=& w - \frac{\hat{\tau}}{\zeta \tau}\Big(w -\mathbb{E}_{\Delta,T,Z_0,Q}\Big[Z_0\xi\Big]\Big) \\ %\theta_0 \frac{\hat{\tau}}{\zeta}\mathbb{E}_{\Delta,T,Z_0,Q}\Big[\dot{g}(\xi, \Lambda(T), \Delta)\dot{g}_0(\theta_0 Z_0, \Lambda_0(T), \Delta)\Big]\\
\label{eq : rs5}
    \zeta \tau /\hat{\tau}&=&\mathbb{E}_{\Delta,T,Z_0,Q}\Big[Q\xi\Big] \\
\label{eq : rs6}
    \zeta \hat{v}^2 &=& \frac{\hat{\tau}^2}{\tau^2}\mathbb{E}_{\Delta,T,Z_0,Q}\Big[(\xi - wZ_0 - vQ)^2\Big]
\end{eqnarray}
where 
\begin{eqnarray}
\label{def : xi}
    \fl \xi&:=& {\rm prox}_{g(.,\Lambda(T),\Delta)}(wZ_0 + vQ, \tau)=  w Z_0 + v Q + \tau \Delta - W_0\Big(\tau\rme^{\tau\Delta + wZ_0 + v Q}\Lambda(T)\Big) \\
\label{def : varphi}
    \fl \varphi &:=&  {\rm prox}_{r(.)}(\hat{w}\beta_0/\theta_0 + \hat{v}Z, \hat{\tau})\ ,
\end{eqnarray}
and $Z\sim \mathcal{N}(0,1)$.

\subsection{Interpretation}
We show in \ref{app_sec : distributions} that 
\begin{eqnarray}
      &&\mathcal{P}(\Delta, t, h):= \lim_{n\rightarrow \infty }\mathbb{E}_{\mathcal{D}}\Big[\frac{1}{n}\sum_{i=1}^n \delta(t - T_i)\delta_{\Delta, \Delta_i}\delta (x - \mathbf{X}_i'\hat{\bbeta})\Big]  \\
     &&\mathcal{P}_{\varphi}(x) := \lim_{n\rightarrow \infty }\mathbb{E}_{\mathcal{D}}\Big[\frac{1}{p}\sum_{k=1}^p \delta (x - \mathbf{e}_{k}'\hat{\bbeta})\Big]  \ .
\end{eqnarray}
admit the following expressions in the Replica Symmetric ansatz
\begin{eqnarray}
\label{def : law_xi_rs}
    \mathcal{P}_{\xi}(\Delta, t, h) &=& \mathbb{E}_{\Delta', T', Z'_0, Q'} \Big[\delta(t - T')\delta_{\Delta, \Delta'}\delta\big(h - \xi_{\star}\big)\Big]\\
\label{def : law_varphi_rs}
    \mathcal{P}_{\varphi}(x) &=& \mathbb{E}_{\beta_0, Z} \Big[\delta\big(x - \varphi_{\star}\big)\Big] 
\end{eqnarray}
where $\xi_{\star}, \varphi_{\star}$ are as in (\ref{def : xi},\ref{def : varphi}), computed at the fixed point of the RS equations, i.e. at $w_{\star}, v_{\star}, \tau_{\star}, \hat{w}_{\star}, \hat{v}_{\star}, \hat{\tau}_{\star}$ and $\Lambda_{\star}(.)$ is the function solving (\ref{eq: self_const_Lambda_1}, \ref{eq: self_const_Lambda_2}) evaluated at the fixed point.
Hence, in the proportional regime and for $\mathbf{X}_i \sim \mathcal{N}(\bm{0}, \frac{1}{p}\bm{I}_p)$, we have
\begin{eqnarray}
\label{def : av_equivalence_xi}
    &&\lim_{n\rightarrow \infty } \mathbb{E}_{\mathcal{D}}\Big[\frac{1}{n} \sum_{i=1}^n \ell_1(\mathbf{X}_i'\hat{\bbeta}_n, \hat{\Lambda}_n(T_i), \Delta_i, T_i)\Big] = \nonumber\\
    &&\hspace{2cm }\sum_{\Delta \in \{1,0\}}\int \mathcal{P}_{\xi}(\Delta, t, h) \ell_1 (h, \Lambda_{\star}(t), \Delta, t) \rmd t  \rmd h, \\
\label{def : av_equivalence_varphi}
   &&\lim_{p\rightarrow \infty } \mathbb{E}_{\mathcal{D}}\Big[\frac{1}{p}\sum_{\mu=1}^p \ell_2(\mathbf{e}_{\mu}'\hat{\bbeta}_n)\Big] =  \int \mathcal{P}_{\varphi}(x) \ell_2(x)\ \rmd x \ , 
\end{eqnarray}
for any \say{reasonable} functions $\ell_1, \ell_2$. This can be compactly re-written as 
\begin{eqnarray}
\label{def : equivalence_xi}
    \bxi:= {\rm prox}_{g(.,\Lambda(\mathbf{T}),\Delta)}(w\mathbf{Z}_0 + v\mathbf{Q}, \tau)&\underset{n\rightarrow\infty}{\overset{d}{\approx}}& \bX\hat{\bbeta}_n\\
\label{def : equivalence_varphi}
    \bvarphi:= {\rm prox}_{r(.)}(\hat{w}\bbeta_0/\theta_0 + \hat{v}\mathbf{Z}, \hat{\tau}) &\underset{n\rightarrow\infty}{\overset{d}{\approx}}& \hat{\bbeta}_n \ ,
\end{eqnarray}
where $\underset{n\rightarrow\infty}{\overset{d}{\approx}}$ means approximately equal in distribution in the limit and  $\mathbf{Z}_0 = (Z_{0,1}, \dots, Z_{0,n})$, $\mathbf{Q} = (Q_1, \dots, Q_n)$ , $\mathbf{Z} = (Z_1, \dots, Z_n)$ are indipendent standard Gaussian vector and the proximal operators act component-wise. 
In turn (\ref{def : av_equivalence_varphi}) or (\ref{def : equivalence_varphi}) provides the following interpretation of the values $w_{\star},v_{\star}$ 
\begin{equation}
\label{approx_conv_beta}
   w_n := \frac{1}{\sqrt{p}}\frac{\bbeta_0'\hat{\bbeta}_n}{\sqrt{\bbeta_0'\bbeta_0}} \underset{n\rightarrow \infty}{\approx} w_{\star}, \quad v_n^2 := \frac{1}{p}\hat{\bbeta}_n'\hat{\bbeta}_n - w_n^2 \underset{n\rightarrow \infty}{\approx} v^2_{\star} \ . 
\end{equation}
In words, the random quantity $w_n$ is the overlap between the estimator 
$\hat{\bbeta}_n$ and the truth $\bbeta_0$, whilst the random variable $v_n$ measures the modulus of the component of $\hat{\bbeta}_n$ that is orthogonal to $\bbeta_0$, i.e. the noise in the estimation of $\bbeta_0$. The deterministic values $w_{\star}, v_{\star}$ are the respective asymptotic values.
Similarly (\ref{def : av_equivalence_xi}) or (\ref{def : equivalence_xi}) implies that the Replica Symmetric functional order parameter $\Lambda(.)$ is equivalent to the Nelson - Aalen estimator computed over a large (ideally infinite) data-set, i.e. 
\begin{equation}
\label{approx_conv_Lambda}
    \hat{\Lambda}_n(.) \underset{n\rightarrow \infty}{\approx} \Lambda_{\star}(.)  \ .
\end{equation}
The identifications (\ref{approx_conv_beta}, \ref{approx_conv_Lambda}), were already present in our recent investigation of the Maximum Partial Likelihood estimator \cite{Massa_24}. However, here they are derived in a more general setting.

The physical meaning of $\hat{w}_{\star}$ and $\hat{v}_{\star}$ can be deduced by Cox-AMP, see Algorithm \ref{alg : amp_cox}.
Since the large $n$ behaviour of Cox-AMP is predicted by the RS equations, we see that, at the fixed point, the local field 
\begin{equation}
\label{def : local_field}
\hspace{-1cm}\bm{\psi}_n :=  \hat{\bbeta} - \hat{\tau}\mathbf{X}'\dot{\mathcal{M}}_{g(. , \hat{\Lambda}_n(\bT), \bm{\Delta})}(\bm{\xi}, \tau) =  \hat{\bbeta} - \hat{\tau}\mathbf{X}'\dot{g}(\mathbf{X}\hat{\bbeta}_n, \hat{\Lambda}_n(\bT), \bm{\Delta})\ , 
\end{equation}
has a Gaussian distribution
\begin{equation}
\label{local_field_distribution}
    \bm{\psi}_{n} \overset{d}{\approx} \hat{w}_{\star} \bbeta_0/\theta_0 + \hat{v} \mathbf{Z}, \quad \mathbf{Z} \sim \mathcal{N}(\bm{0}_p, \bm{I}_p) \ .
\end{equation}
This, in turn, implies that 
\begin{equation}
   \hspace{-1cm} \hat{w}_n := \frac{1}{\sqrt{p}}\frac{\bbeta_0'\bm{\psi}_{n}}{\sqrt{\bbeta_0'\bbeta_0}}\underset{n\rightarrow \infty}{\approx} \hat{w}_{\star}, \quad \hat{v}_n^2 := \frac{1}{p}\bm{\psi}_{n}'\bm{\psi}_{n} - \hat{w}_n^2 \underset{n\rightarrow \infty}{\approx} \hat{v}^2_{\star} \ . 
\end{equation}

\subsection{Application: approximate cross validation}
\label{subsec : application}
By definition of proximal operator (\ref{def : proximal}),
\begin{equation}
    \bxi = w \mathbf{Z}_0 + v \mathbf{Q} - \tau g(\bxi, \mathbf{T}) \ ,
\end{equation}
hence the identification (\ref{def : equivalence_xi}) suggests that 
\begin{equation}
\label{def : tilde_bxi}
    \tilde{\bxi} := \mathbf{X}\hat{\bbeta}_n + \tau_{\star}g(\mathbf{X}\hat{\bbeta}_n, \mathbf{T})\overset{d}{\approx}  w \mathbf{Z}_0 + v \mathbf{Q} \ .
\end{equation}
If one is able to estimate $\tau_{\star}$, and we can, as we shall soon show, then the expression above is useful to estimate the generalization error. To see this, notice that the linear predictor for a fresh observation $\mathbf{X}_{\star}\sim \mathcal{N}(\bm{0}, \bm{I}_p / p)$ is 
\begin{equation}
    \mathbf{X}_{\star}'\hat{\bbeta}_n \overset{d}{=} w_n Z_{0, \star} + v_n Q_{\star} \overset{d}{\underset{n\rightarrow\infty}{\longrightarrow}} w_{\star} Z_{0, \star} + v_{\star} Q_{\star}, \  
\end{equation}
with $Z_{0, \star}, Q_{\star}$ independent standard Gaussian random variables, $\overset{d}{=}$ indicates equality in distribution and $\overset{d}{\underset{n\rightarrow\infty}{\longrightarrow}}$ indicates convergence in distribution. 
The statement above is rigorous if $w_n \overset{P}{\underset{n\rightarrow\infty}{\longrightarrow}} w_{\star}$ and $v_n \overset{P}{\underset{n\rightarrow\infty}{\longrightarrow}} v_{\star}$, because of Slutsky's lemma \cite{Van_der_Vaart_2000}. 

This implies that each  component of $\tilde{\bxi}$ has the same distribution as $\mathbf{X}_{\star}'\hat{\bbeta}_n$. Thus, one can evaluate generalization metrics by using the training responses $\bT, \bm{\Delta}$ with linear predictors $\tilde{\bxi}$. For instance, the Test Harrel's C-index can be estimated as follows 
\begin{equation}
    \label{def : CV_C_index}
        \widehat{HC}_{RSCV} :=\frac{ \sum_{i=1}^{n} \Delta_i \sum_{j=1}^{n} \Theta(T_j-T_i) \Theta\Big(\tilde{\xi}_j-\tilde{\xi}_i\Big)}{ \sum_{i=1}^{n_{test}} \Delta_i \sum_{j=1}^{n_{test}} \Theta(T_j-T_i) }  \ ,
\end{equation}
where $\tilde{\bxi}$ is defined in (\ref{def : tilde_bxi}), and $\tilde{\xi}_i$ is the i-th component of $\tilde{\bxi}$.
The expression on the right-hand side above resembles a leave one cross validated estimate of the C -index, where the leave one out linear predictors are replaced by $\tilde{\bxi}_i$. Thus, we dub the quantity above \say{RSCV} for Replica Symmetric Cross Validated. An application of this idea to simulated data can be seen in figure (\ref{fig : elbow_plot}).

\section{Numerical solution of the RS equations with known data-generating process}
\label{sec : simulations}
In this section we compare the solution of the RS equations (\ref{eq: self_const_Lambda_1}, \ref{eq: self_const_Lambda_2}, \ref{eq : rs1}, \ref{eq : rs2}, \ref{eq : rs3}, \ref{eq : rs4}, \ref{eq : rs5}, \ref{eq : rs6}) with the finite sample size metrics computed from synthetic data when the data generating process is perfectly known. 
This is done to benchmark the prediction of the Replica Symmetric theory. 

We first detail on the simulation protocol. 
In all the simulations of this section, the true associations are generated as follows
\begin{equation}
    \bbeta_0 = \Big(\theta_0\  \sqrt{p}\  \bm{U}_{s}, \bm{0}_{p-s}\Big), \quad \bm{U}_s \sim {\rm Unif}(\mathbb{S}_{s-1}) \ . 
\end{equation}
When $p,s\rightarrow\infty$ we have that 
\begin{equation}
    \fl \mathcal{P}_{\beta_0}(x) = \nu \frac{1}{\sqrt{2\pi}\sigma}\rme^{-\frac{1}{2} x^2/\sigma^2 } + (1-\nu) \delta(x), \quad \nu = s/p \in [0,1] , \quad \sigma = \theta_0/\sqrt{\nu} \ .
\end{equation}
where $\mathcal{P}_{\beta_0}$ is the limiting empirical distribution of the entries of $\bbeta_0$, defined in (\ref{def : p_beta0}).
The sparsity of the signal is controlled by $\nu$ which is the fraction of \say{active} covariates that is actually correlated with the outcome. 
Under this assumption, and with the elastic net regularization (\ref{def : elastic_net}), the replica symmetric equations can be simplified further, see \ref{app_sec : computation_enet}, to a somewhat simpler form
\begin{eqnarray}
    w &=& 2\frac{1}{1+\eta\hat{\tau}}\hat{w} \Phi(\chi_1)\\
    \tau &=& 2\frac{1}{1+\eta\hat{\tau}}\hat{\tau}\Big\{ \nu \Phi(\chi_1) + (1-\nu) \Phi(\chi_0)\Big\}\\
    \frac{1}{2}(v^2 + w^2) &=&  \frac{1}{(1+\eta\hat{\tau})^2}\Big\{ \nu \Big((1 + 1/\chi_1^2) \Phi(\chi_1) - \phi(\chi_1)\Big) +\nonumber \\
    &+& (1-\nu) \Big( (1 + 1/\chi_0^2) \Phi(\chi_0) - \phi(\chi_0)\Big)\Big\}\\
     \hat{w} &=& w - \frac{\hat{\tau}}{\zeta \tau}\Big(w -\mathbb{E}_{\Delta,T,Z_0,Q}\Big[Z_0\xi\Big]\Big)\\
    v(1 -\zeta\tau/\hat{\tau}) &=&\mathbb{E}_{\Delta,T,Z_0,Q}\Big[Q\xi\Big] \\
    \zeta \hat{v}^2 &=& \frac{\hat{\tau}^2}{\tau^2}\mathbb{E}_{\Delta,T,Z_0,Q}\Big[(\xi - wZ_0 - vQ)^2\Big]
\end{eqnarray}
where $\xi$ is defined in (\ref{def : xi}), 
\begin{equation}
    \chi_0 := \frac{\alpha\hat{\tau}}{\hat{v}}, \quad \chi_1 := \frac{\alpha\hat{\tau}}{\sqrt{\hat{v}^2 + \hat{w}^2/\nu}}, 
\end{equation}
and we indicated with $\phi(x) := \exp\{-\frac{1}{2}x^2\}/\sqrt{2\pi}$ the density function and with $\Phi(x) := \int_x^{+\infty} \phi(t) \rmd t$ the complementary distribution function of a standard normal random variable.
The numerical solution of these self -consistent equations is obtained via fixed point iteration. The integrals are approximated as averages of populations of size $5 \cdot 10^3$.
The latent survival data are generated from a Log-logistic proportional hazard model
\begin{eqnarray}
     Y|\mathbf{X} &\sim&   -\frac{\rmd }{\rmd  t }  S_0(t|\mathbf{X}), \nonumber \\
    S_0(t|\mathbf{X}) &=& \exp\{- \Lambda_0(t)\exp(\mathbf{X}'\bbeta_0)\}, \ \Lambda_0(t) := \log \Big(1 + \rme^{\phi}_0 t^{\rho_0}\Big) \ , 
\end{eqnarray}
and the latent censoring time is sampled uniformly in the interval $[\tau_1, \tau_2]$, i.e.
\begin{equation}
    C \sim {\rm Unif}[\tau_1, \tau_2]\ . 
\end{equation}
The observations $\Delta, T$ are then generated according to 
\begin{equation}
    \Delta = \bm{1}[Y<C], \quad Y = \min\{Y,C\} \ .
\end{equation}
To make the setting above concrete, we choose: $\theta_0 = 1, \phi_0 = -\log 2, \rho_0 = 2$ and $\tau_1 = 1.0, \tau_2 = 2.0$, with $\rho = 0.75$ and $\nu = 0.005$.
The quantities computed from simulations are always displayed as dots, indicating the averages over $20$ repetitions (optimization with COX-AMP), with error bars, the corresponding finite sample standard deviations. 

\begin{figure}[H]
\centering
\includegraphics[scale = 0.9]{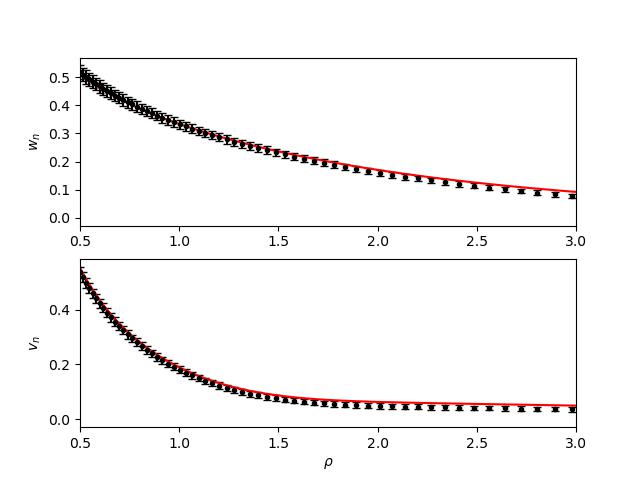}
\caption{\textbf{Comparison} between the numerical solution of the RS equations against their finite sample counter part $w_n := (\hat{\bbeta}_n'\bbeta_0/\|\bbeta_0\|)/\sqrt{p}$, $
v_n := \sqrt{\|\hat{\bbeta}_n^2/p - w_n^2\|}$. Notice that here we use the knowledge of $\bbeta_0$ to compute $w_n, v_n$. Here $\zeta = 2.0, \nu = 0.005$ and $p = 2000$, i.e. $n = 1000$, $s =10$ .} 
\label{fig : amp}
\end{figure}

We have seen that the prediction of the RS theory are indeed correct. However, the RS equations cannot be used in practice, since they require perfect knowledge of the data-generating process. To overcome this limitation we need to relate the order parameters to \emph{practically} observable quantities, i.e. quantities that can be computed without knowing the data generating process, i.e. $\bbeta_0 $ and $\Lambda_0(.)$. We show that this is indeed possible, and in a quite straightforward manner, when the RMPLE is obtained via the Cox - AMP Algorithm \ref{alg : amp_cox}.

\section{Inferring the RS order parameters via Cox-AMP without solving the RS equations }
\label{sec: method_amp}
The COX-AMP algorithm returns $\tau_n$ and $\hat{\tau}_n$, which we interpret as finite sample size realizations of $\tau_{\star}$ and $\hat{\tau}_{\star}$. 
This means that $\hat{v}$ can be estimated as $v_n$ from the data, by approximating 
\begin{equation}
     \zeta \hat{v}^2 = \hat{\tau}^2 \mathbb{E}_{T,Z_0,Q}\Big[\dot{g}\big(\xi, \Lambda(T), \Delta\big)^2\Big] \ , 
\end{equation}
with 
\begin{equation}
     \hat{v}_n^2 = \big\langle \ddot{g}(\mathbf{X}\hat{\bbeta}_n, \hat{\Lambda}_n, \bT)\big\rangle \hat{\tau}_n^2  / \zeta \ ,
\end{equation}
via (\ref{def : equivalence_xi}), since $\zeta$ is known.
Above we have indicated with $\langle\rangle$ the empirical average operator, i.e. $\langle \bm{a}\rangle = \sum_{k=1}^d a_k / d$ for a vector $\bm{a}\in \mathbb{R}^d$. The functions $g,\dot{g}, \ddot{g}, \hat{\Lambda}$ are taken to act component wise.

\begin{figure}[H]
\centering
\includegraphics[width = \textwidth]{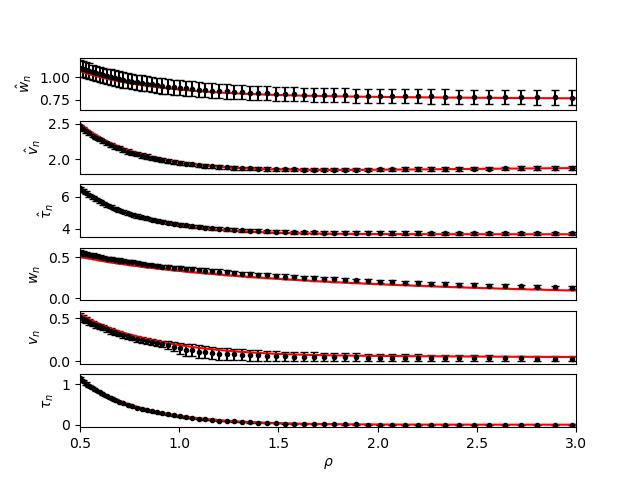}
\caption{\textbf{Cox-AMP with elastic net regularization} The estimators of $\hat{w}_{\star}, \hat{v}_{\star}, \hat{\tau}_{\star}, w_{\star}, v_{\star}, \tau_{\star}$  summarized by a black errorbar plot against the true values in red, along a regularization path (different values of $\alpha$). The data are generated from a Log-logistic proportional hazards model $\Lambda_0(t) = \log(1+t^2/2)$ with uniform censoring between $\tau_1 = 1.0$ and $\tau_2 = 2.0$. The associations  $\bbeta_0$ are sampled as from a Gauss-Bernoulli distribution $\mathcal{P}_{\beta_0}(x) = (1-\nu) \delta(x) + \nu \exp\{-x^2/2\sigma^2\} / \sqrt{2\pi}\sigma$, with sparsity $\nu = 0.005$ and $\sigma$ adjusted to have $\|\bbeta_0\|^2_2/p = \theta^2_0 = 1.0$. The L1 ratio is $\lambda = 0.75$. } 
\label{fig : amp_order_parameters}
\end{figure}

We have already noticed in the Interpretation subsection that the local field $\bm{\psi}_n$ of the Cox-AMP algorithm, defined in (\ref{def : local_field}), has a Gaussian distribution (\ref{local_field_distribution}). By matching the variance, this implies that
\begin{equation}
\label{def : amp_equation}
    \frac{1}{p}\big\|\hat{\bbeta} - \hat{\tau} \mathbf{X}'\dot{\mathcal{M}}_{g(. , \hat{\Lambda}_n(\bT), \bm{\Delta})}(\mathbf{X}\hat{\bbeta}_n, \tau_{\star}) \big\|^2_2 =  \hat{w}_n^2 + \
    \hat{v}_n^2 \ . 
\end{equation}
From which we infer $\hat{w}_n$. 
We now address the slightly more complicated issue of estimating $w_{\star}$ and $v_{\star}$ from the data. 
It can be seen that, by definition of proximal operator for $\xi$, we have 
\begin{eqnarray}
    \fl \zeta \hat{v}^2 \tau^2/ \hat{\tau}^2 &=& \tau^2\mathbb{E}_{\Delta,T,Z_0,Q}\Big[\dot{g}(\xi,\Lambda(T), \Delta)^2\Big]  =\mathbb{E}_{\Delta,T,Z_0,Q}\Big[\big(\xi - w Z_0 - v Q\big)^2\Big] =  \nonumber \\
    \fl &=& \mathbb{E}_{\Delta,T,Z_0,Q}\Big[\xi^2\Big] - 2 w \mathbb{E}_{\Delta,T,Z_0,Q}\Big[Z_0 \xi\Big]  - 2 v \mathbb{E}_{\Delta,T,Z_0,Q}\Big[Q \xi\Big] + w^2 + v^2= \nonumber\\
    \fl &=& \mathbb{E}_{\Delta,T,Z_0,Q}\Big[\xi^2\Big] + 2 w \zeta \frac{\tau}{\hat{\tau}}(w - \hat{w})  + 2 v^2 \zeta \frac{\tau}{\hat{\tau}}- w^2 - v^2 =\nonumber \\
    \fl &=& \mathbb{E}_{\Delta,T,Z_0,Q}\Big[\xi^2\Big] - (w^2 + v^2) (1 - 2 \zeta \frac{\tau}{\hat{\tau}}) - 2 w \hat{w} \zeta \frac{\tau}{\hat{\tau}}  \ .
\end{eqnarray}
The above implies 
\begin{equation}
    \mathbb{E}_{\Delta,T,Z_0,Q}\Big[\xi^2\Big] = \zeta \hat{v}^2 \tau^2/ \hat{\tau}^2 + (w^2 + v^2) (1 - 2 \zeta \frac{\tau}{\hat{\tau}}) +  2 w \hat{w} \zeta \frac{\tau}{\hat{\tau}}  \ . 
\end{equation}
On the other hand, the relationship 
\begin{equation}
    \mathbb{E}_{\Delta,T,Z_0,Q}\Big[\big(\xi + \tau \dot{g}(\xi, \Lambda(T), \Delta)\big)^2\Big] = w^2 + v^2  \ , 
\end{equation}
can be deduced by the definition of proximal operator. 
Hence 
\begin{equation}
    w \hat{w} \zeta \frac{\tau}{\hat{\tau}} = \frac{1}{2} \mathbb{E}_{\Delta,T,Z_0,Q}\Big[\xi^2\Big]  - \frac{1}{2} \zeta \hat{v}^2 \tau^2 / \hat{\tau}^2 - \frac{1}{2} (w^2 + v^2) (1 - 2 \zeta \tau / \hat{\tau})  \ .
\end{equation}
The finite sample counterpart of the previous equations read
\begin{eqnarray}
     \hspace{-0.5cm}&&\frac{1}{n}\|\bX\hat{\bbeta}_n + \tau_n \dot{g}(\bX\hat{\bbeta}_n, \hat{\Lambda}_n(\bT), \bm{\Delta})\|_2^2 = w_n^2 + v_n^2 \\
    \hspace{-0.5cm}&&w_n \hat{w}_n \zeta \frac{\tau_n}{\hat{\tau}_n} = \frac{1}{2n} \|\bX \hat{\bbeta}_n\|_2^2   - \frac{1}{2} \zeta \hat{v}_n^2 \tau_n^2 / \hat{\tau}_n^2 - \frac{1}{2} (w_n^2 + v_n^2) (1 - 2 \zeta \tau_n / \hat{\tau}_n)  \ .
\end{eqnarray}
From which one can solve for $v_n$ and $w_n$. The remarkable performance of the estimators obtained as described above in estimating the replica symmetric order parameters $w_{\star}, v_{\star}, \tau_{\star}, \hat{w}_{\star}, \hat{v}_{\star}, \hat{\tau}_{\star}$ can be appreciated in figure(\ref{fig : amp_order_parameters}). There we show the estimators summarized as errorbar plots, computed over $20$ realizations,  when $\zeta = 2.0$ and $p = 2000$, for the elastic net regularization (\ref{def : elastic_net}) when the L1 ratio is $\lambda = 0.75$. 

We have shown that the RS order parameter can be easily estimated from the estimators $\hat{\bbeta}_n$ and $\hat{\Lambda}_n$ as computed from the COX-AMP algorithm. The key property of the AMP algorithm is that it estimates also $\hat{\tau}$ and $\tau$ (which are not available in other optimization methods), which in turn allows to estimate the local field $\bm{\psi}_{\star}$ from the data. It seems natural to wonder if a similar construction is available when $\hat{\bbeta}_n$ is computed via another method, e.g. Coordinate-wise Descent. We show in the next section how this can be done.

\section{Inferring the RS order parameters via Cox-CD without solving the RS equations }
\label{sec: method_cd}
\begin{figure}[H]
\centering
\includegraphics[width = \textwidth]{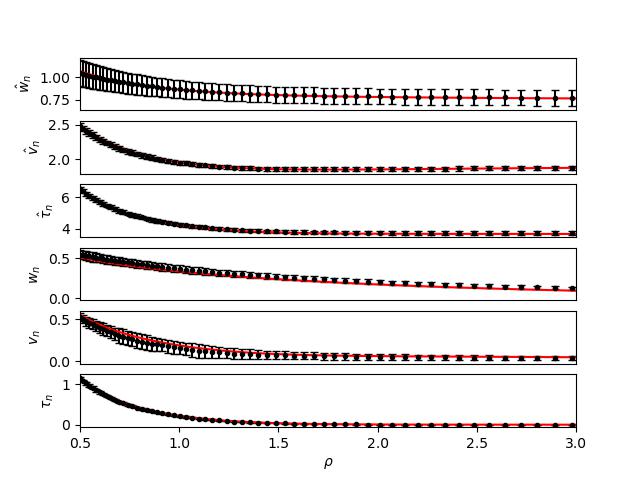}
\caption{\textbf{Cox-CD with elastic net regularization} The estimators of $\hat{w}_{\star}, \hat{v}_{\star}, \hat{\tau}_{\star}, w_{\star}, v_{\star}, \tau_{\star}$  summarized by a black errorbar plot against the true values in red, along a regularization path (different values of $\alpha$). The data are generated from a Log-logistic proportional hazards model $\Lambda_0(t) = \log(1+t^2/2)$ with uniform censoring between $\tau_1 = 1.0$ and $\tau_2 = 2.0$. The associations  $\bbeta_0$ are sampled as from a Gauss-Bernoulli distribution $\mathcal{P}_{\beta_0}(x) = (1-\nu) \delta(x) + \nu \exp\{-x^2/2\sigma^2\} / \sqrt{2\pi}\sigma$, with sparsity $\nu = 0.005$ and $\sigma$ adjusted to have $\|\bbeta_0\|^2_2/p = \theta^2_0 = 1.0$. The L1 ratio is $\lambda = 0.75$. } 
\label{fig : cd}
\end{figure}

Within the Coordinate-wise Descent algorithm, we do not infer the value of $\tau$, nor the value of $\hat{\tau}$. Hence we cannot directly observe $\hat{v}_n$ and solve the chain of equations that we have derived in the previous section. However, in some cases it seems possible to infer $\hat{\tau}$ and $\tau$ nevertheless. 
Noticing that 
\begin{eqnarray}
    \hat{v}\frac{\tau}{\hat{\tau}} &=&\mathbb{E}_{Z, \beta_0}\Big[Z \varphi\Big] =\nonumber \\
     &=&\hat{v}\mathbb{E}_{Z, \beta_0}\Big[\prox_{\rmr(.)}'\big(\hat{w}\beta_0 / \theta_0 + \hat{v} Z, \hat{\tau}\big)\Big] \ , 
\end{eqnarray}
where we used integration by parts, also known as Stein's lemma, we obtain an alternative form of (\ref{eq : rs2})
\begin{equation}
    \frac{\tau}{\hat{\tau}} = \mathbb{E}_{Z, \beta_0}\Big[\prox_{\rmr(.)}'\big(\hat{w}\beta_0 / \theta_0 + \hat{v} Z, \hat{\tau}\big)\Big] \ .
\end{equation}
For the elastic net regression for instance, we have 
\begin{equation}
    {\rm prox}_{\rmr(.)}\big(\hat{w} \beta_0/\theta_0 + \hat{v} Z, \hat{\tau}\big) = \frac{1}{1 + \eta \hat{\tau}}{\rm st}\big(\hat{w} \beta_0/\theta_0 + \hat{v} Z, \alpha \hat{\tau}\big)
\end{equation}
and
\begin{eqnarray}
        {\rm prox}_{\rmr(.)}'\big(\hat{w} \beta_0/\theta_0 + \hat{v} Z, \hat{\tau}\big) &=& \frac{1}{1 + \eta \hat{\tau}}{\rm st}'\big(\hat{w} \beta_0/\theta_0 + \hat{v} Z, \alpha \hat{\tau}\big) \nonumber \\
        &=& \frac{1}{1 + \eta \hat{\tau}} \Theta(|\varphi|) \ , 
\end{eqnarray}
where in the last expression we used the fact that 
\begin{equation}
    {\rm st}'\big(\hat{w} \beta_0/\theta_0 + \hat{v} Z, \alpha \hat{\tau}\big) = \Theta\big(|\hat{w} \beta_0/\theta_0 + \hat{v} Z| - \alpha \hat{\tau}\big) = \Theta(|\varphi|) \ .
\end{equation}
So 
\begin{equation}
    \frac{\tau}{\hat{\tau}} = \frac{1}{1 + \eta \hat{\tau}} \mathbb{E}_{Z, \beta_0}\Big[\Theta(|\varphi|)\Big] \ ,  
\end{equation}
which might be solved simultaneously with (\ref{eq : rs5}) for $\hat{\tau}, \tau$. Their finite sample counter-parts read
\begin{eqnarray}
\label{eq : elastic_net_1}
    \frac{\tau_n}{\hat{\tau}_n} &=& \frac{1}{1 + \eta \hat{\tau}_{n}} \|\hat{\bbeta}_n\|_0\\
\label{eq : elastic_net_2}
    \zeta \tau_n / \hat{\tau}_n &=& \Big\langle \frac{\tau_n\ddot{g}(\bX\hat{\bbeta}_n, \hat{\Lambda}_n(\bT), \bm{\Delta})}{ 1 + \tau_n \ddot{g}(\bX\hat{\bbeta}_n, \hat{\Lambda}_n(\bT), \bm{\Delta})}\Big\rangle \ , 
\end{eqnarray}
where we indicated with $\|\mathbf{x}\|_0$ the zero \say{norm} of $\mathbf{x}$, i.e. the number of non-null entries of $\mathbf{x}$.
The system above can be easily solved by first solving (\ref{eq : elastic_net_1}) for $\hat{\tau}_n$ as a function of $\tau_n$, i.e. 
\begin{equation}
    \hat{\tau}_n  = \frac{\tau_n}{ \|\hat{\bbeta}_n\|_0 - \eta \tau_n}  \, 
\end{equation}
which gives by substitution
\begin{equation}
    \zeta \big(\|\hat{\bbeta}_n\|_0 - \eta \tau_n\big)  = \Big\langle \frac{\tau_n\ddot{g}(\bX\hat{\bbeta}_n, \hat{\Lambda}_n(\bT), \bm{\Delta})}{ 1 + \tau_n \ddot{g}(\bX\hat{\bbeta}_n, \hat{\Lambda}_n(\bT), \bm{\Delta})}\Big\rangle \ , 
\end{equation}
which can be solved numerically for $\tau_n$, by Newton method for instance.
We see in figure(\ref{fig : cd}) that the resulting estimators, summarized as black error bar plots, are close to the quantity to be estimated (i.e. the fixed point solution of the RS equations in red).  It is indeed difficult to see differences between the estimators obtained via the CD algorithm (\ref{fig : cd}) and the ones obtained via the AMP algorithm (\ref{fig : amp_order_parameters}), as it should. 
We conclude that the Cox-AMP algorithm is not the unique mean to estimate the values of the RS order parameters, i.e. it is not strictly needed. However, the additional equation (\ref{def : amp_equation}), which is necessary to close the estimating equations, is derived via the insight provided by the AMP algorithm and its statistical physics interpretation. 

\section{Conclusion}
Under the assumption of centered Gaussian covariates with variance $1/p$ and convex regularization function, we showed via simulations that the Replica Symmetric theory is capable of correctly predicting the behaviour of the Regularized Maximum Partial Likelihood Estimator in the proportional regime when the number of covariates $p$, the number of samples $n$ and the number of active covariates $s$ diverge proportionally. When the data generating model (i.e. true associations, $\bbeta_0$, and true base hazard rate $\lambda_0$) is perfectly known, the theory gives the values of several observables of interest, such as the Bias and the Mean Squared Error, as a function of six order parameters $w_{\star}, v_{\star}, \tau_{\star}, \hat{w}_{\star}, \hat{v}_{\star}, \hat{\tau}_{\star}$ and $\theta_0 := \|\bm{\Sigma}_0^{1/2}\bbeta_0\|_2$. However, the data generating process is not available in applications. We show here that one can estimate the RS order parameters \emph{solely} from the data, \emph{without} the knowledge of the data generating process. \\
We proposed an AMP algorithm tailored to the Cox model in section \ref{sec : AMP_CD} and note that this suggests the construction of a local field
\begin{equation}
\bm{\psi}_{\star} :=  \hat{\bbeta}_n - \hat{\tau}_{\star} \mathbf{X}'\dot{\mathcal{M}}_{g(. , \hat{\Lambda}_n(\bT), \bm{\Delta})}(\bm{\xi}, \tau_{\star})  \ , 
\end{equation}
which is predicted to have a Gaussian distribution by comparison with the RS equations, i.e.
\begin{equation}
    \bm{\psi}_{\star} \overset{d}{=} \hat{w}_{\star} \bbeta_0/\theta_0 + \hat{v}_{\star} \mathbf{Z}, \quad \mathbf{Z} \sim \mathcal{N}(\bm{0}_p, \bm{I}_p) \ .
\end{equation}
This identification allows estimating all the RS order parameters directly from the data, when the RMPLE $\hat{\bbeta}_n$ is obtained from the COX-AMP algorithm since in that case $\hat{\tau}_{\star}$ is observed, as explained in sections \ref{sec: method_amp}. Afterwards, we noticed that, in some cases, it is possible to estimate the local field $\bm{\psi}_{\star}$ also from the RMPLE obtained by any other method, e.g. Coordinate-wise Descent, i.e. without the need of explicitly using the AMP algorithm, in section \ref{sec: method_cd}.

When the regularization is non convex, Replica Symmetry Breaking might occur, see for instance \cite{Sakata_2023} where the author derived the 1-RSB free energy for GLMs. In this case, additional parameters needs to be estimated and the theory, in its actual form, is not be observable, i.e. one can not compute the values of the RS order parameters solely from the data. The extension of our results to this more challenging setting will be the subject of future investigations.

The COX-AMP algorithm is derived, as the original AMP \cite{Donoho_2009} and GAMP \cite{Rangan_2010}, under the assumption of i.i.d. centered Gaussian covariates with variance $1/p$. Because of universality \cite{Bayati_2015}, we expect the algorithm and proposed methodology, to \say{work} when the (centered, normalized and eventually rescaled) covariates are uncorrelated, sub-Gaussian random variables (see \cite{Vershynin_2018, Boucheron_2013} for a definition of sub-Gaussian random variable). The same is true for the Replica derivation, which is conjectured (rightly so, as we show via simulations) to describe the behaviour of the minimizer of the Regularized Partial Likelihood objective (\ref{def : ppl}).
Although we focused on this idealized setting, which limits the applicability of the theory presented here, the present study achieves the goal of estimating the RS order parameters solely from the data, by-passing the need of solving the RS equations, and paves the road to future applications with real data. This might include, but is not restricted to, the tuning of the regularization parameters, for instance by maximizing the ratio of signal to noise (as quantified by the RS order parameters $w_{\star}$ and $v_{\star}$ respectively) or the computation of the prediction error.

However, further work is still required before applications can materialize.
First, it is desirable to account for correlations between the covariates. If the covariates are correlated, the RS order parameters depend implicitly on the true covariance matrix, and it is not yet clear how one can estimate the order parameters solely from the data, without the knowledge of the true covariance matrix. Furthermore, the COX-AMP algorithm is not valid, or at least is not expected to converge. 
%For this reason, it is not yet clear how one can estimate the order parameters solely from the data in presence of correlation between the covariates. 
There exist several generalizations of GAMP that are provably more stable (but also more complex), a relatively recent example is  Vector Approximate Message Passing (VAMP) \cite{Rangan_2018, Schniter_2016}. It is of future interest to study similar generalizations of COX-AMP, which might guarantee wider applicability of the algorithm.  
A second desirable requirement is to remove the strong requirement of Gaussian covariates. 
Based on the works \cite{Adomaityte_2023, Adomaityte_2024}, the Replica calculation might be extended to incorporate elliptical and heavy tailed covariates. Other recent studies explored the case of covariates with a Gaussian mixture distribution \cite{Loureiro_2021} and universality claims within this setting \cite{Pesce_2023}. In both settings mentioned before, we necessarily need to estimate additional parameters (the ones of the  mixture in \cite{Loureiro_2021} or of the super-statistical distribution in \cite{Adomaityte_2024}). It is an interesting problem for the future to determine if it is possible to extend the methodology proposed here to incorporate these more general settings, this would increase the applicability of the present results.  

Whilst the main subject of the manuscript is the Cox regression model, the methodology used here to estimate the RS order parameters of the theory can be easily extended \say{mutatis mutandis} to arbitrary GLMs with eventual additional nuisance parameters. In the GLM case, for instance, the function $g$ would have the role of the minus logarithm of the likelihood and the methodology would even simplify, as we would not have the functional equation for the cumulative hazard. It seems that the only requirements for the methodology to be applicable are: 1) that the conditional distribution of the response depends on a single projection of $\bX_i$, i.e. on $\bX_i'\bbeta$ and 2) the minus log-likelihood of the model is convex in the parameters. This is in agreement with recent results \cite{Bellec_2024}, which are however, obtained via an alternative route.   

The GLM category includes survival models, like the parametric Accelerated Failure Time (AFT) models, which do not fall in the PH class. Other more complex survival analysis models, for instance the competing risks models or Cox models with time varying coefficients, might be, in principle, studied via the statistical mechanics formulation used here. This would complicate the Replica derivation and the related AMP algorithm, in that the parameters to be inferred would mathematically represented as matrices rather than vectors. This requires a non trivial and interesting generalization of the theory presented here, which has been recently obtained in the context of classification in \cite{Cornacchia_2023, Cui_2025}. We consider this to be an interesting topic for future investigations.

\section*{Data availability statement}
The python scripts for generating the data, solving the RS equations, fitting the elastic net Cox model and computing the RS order parameters from the data can be found at  \url{https://github.com/EmanueleMassa/Regularized_Cox_model} . 

\section*{Acknowledgements}
A.C.C.C. is grateful for support from CRUK, under grant EDDCPGM /\ 100001.
The authors are grateful to the three anonymous referees that greatly contributed to improve the clarity and exposition of the manuscript. 

\section*{References}

\bibliographystyle{unsrt}
\bibliography{refs}

\appendix

\section{Replica derivation}
\label{app_sec : replica_computation}
The penalized Cox model is defined by the following energy function
\begin{equation}
    \label{Hamiltonian}
     \mathcal{H}_n(\bbeta|\mathcal{D}) = \sum_{i=1}^n \Delta_i \Big[\log\Big(\frac{1}{n}\sum_{j=1}^n \Theta(T_j - T_i)\rme^{\mathbf{X}_i'\bbeta}\Big) - \mathbf{X}_i'\bbeta \Big] +r(\bbeta) 
\end{equation}
where we indicate with $r$ the elastic net penalization
\begin{equation}
    r(\bbeta) = \frac{1}{2} \eta  \|\bbeta\|^2 + \alpha |\bbeta| \ .
\end{equation}
Equivalently we can regard the energy as a functional of the empirical distribution 
\begin{equation}
\label{def : functional_order_parameter}
    P_n(\Delta,t,h|\bbeta,\mathcal{D}) = \frac{1}{n}\sum_{i=1}^n \delta(t - T_i)\delta_{\Delta,\Delta_i} \delta(h - \mathbf{X}_i' \bbeta) 
\end{equation}
where we indicated with $\mathcal{D}$ the data-set, i.e. the set of couples $\{T_j,\mathbf{X}_j\}_{j=1}^n$, which plays here the role of the disorder. Explicitly 
\begin{eqnarray}
     \hspace{-1cm}\mathcal{H}\Big[P_n(.|\bbeta,\mathcal{D} )\Big] &=& n \mathcal{E}\Big[P_n(.|\bbeta,\mathcal{D} )\Big] + r(\bbeta) ,  \\
     \hspace{-0.5cm}\mathcal{E}\Big[P_n(.|\bbeta,\mathcal{D} )\Big] &=& \sum_{\Delta = \pm 1}\int  \Delta \Big(\log  S_n(t|\bbeta,\mathcal{D}) - h \Big) P_n(\Delta,t,h|\bbeta,\mathcal{D}) \rmd h \rmd t , \\
     \hspace{-1cm} S_n(t|\bbeta,\mathcal{D}) &=& \sum_{\Delta' = \pm 1} \int \Theta(t-t') \rme^{h'} P_{n}(\Delta',t,'h'|\bbeta,\mathcal{D})  \rmd h'\rmd t' \ .
\end{eqnarray}
It is well known that the information content of inference is quantified by the free energy 
\begin{equation}
    \fl f(\gamma) := -\lim_{n\rightarrow \infty} \frac{1}{n\gamma} \mathbb{E}_{\mathcal{D}}\Big[\log Z_n(\gamma,\mathcal{D})\Big], \qquad Z_n(\gamma):= \int \rme^{-\gamma  \mathcal{H}[P_n(.|\bbeta,\mathcal{D})]} \rmd \bbeta
\end{equation}
To compute the average of the logarithm we use the replica trick 
\begin{equation}
    f(\gamma) = \lim_{n\rightarrow \infty} \lim_{r\rightarrow 0} f_n^{(r)}(\gamma), \qquad  f_n^{(r)}(\gamma) := -\frac{1}{nr}\log \mathbb{E}_{\mathcal{D}}\Big[Z_n^r(\gamma,\mathcal{D})\Big] \ .
\end{equation}
We refer to  $f_n^{(r)}(\gamma)$ as the replicated free energy.
We now compute this quantity.

\subsection{The replicated free energy}
In order to compute 
\begin{eqnarray}
    \mathbb{E}_{\mathcal{D}}\Big[Z_n^r(\gamma,\mathcal{D})\Big] &=& \int \rme^{-\gamma  \sum_{\alpha=1}^r \mathcal{H}[P_n(.|\bbeta_{\alpha},\mathcal{D})]} \prod_{\alpha=1}^r \rmd \bbeta_{\alpha} \ , 
\end{eqnarray}
We first write $Z_n(\gamma,\mathcal{D})$ in a more convenient form.
Let us introduce the functional delta measure 
\begin{eqnarray} 
\label{def : functional_delta}
     \fl && \delta\Big[\mathcal{P}(.) - P_n(.|\bbeta_{\alpha},\mathcal{D})\Big] \propto \int \rme^{\rmi n\big\{\sum_{\Delta = \pm 1} \int \hat{\mathcal{P}}(\Delta,t,h)\big[\mathcal{P}(\Delta,t,h) - P_n(\Delta,t,h|\bbeta,\mathcal{D})\big] \rmd t \rmd h \big\} }\mathcal{D}\hat{\mathcal{P}} \nonumber \\ 
    \fl  &&\hspace{1cm}= \int \rme^{\rmi n \big\{\sum_{\Delta = \pm 1} \int \hat{\mathcal{P}}(\Delta,t,h)\mathcal{P}(\Delta,t,h)\rmd t \rmd h \big\}- \rmi \sum_{i=1}^n \hat{\mathcal{P}}(\Delta_i,T_i,\mathbf{X}_i'\bbeta_{\alpha}) }\mathcal{D}\hat{\mathcal{P}}
\end{eqnarray}
obtaining
\begin{eqnarray}
     \fl Z_n(\gamma,\mathcal{D})&=& \int \rme^{ n \Big\{ \rmi\sum_{\Delta = \pm 1} \int \hat{\mathcal{P}}(\Delta,t,h)\mathcal{P}(\Delta,t,h)\rmd t \rmd h - \gamma \mathcal{E}\big[\mathcal{P} (.)\big] \Big\} }\times \nonumber \\
     \fl &\times& \bigg\{\int \rme^{- \rmi \sum_{i=1}^n \hat{\mathcal{P}}(\Delta_i,T_i,\mathbf{X}_i'\bbeta) -\gamma r(\bbeta/\sqrt{p})} \rmd \bbeta\bigg\}\ \mathcal{D}\hat{\mathcal{P}} \ \mathcal{D}\mathcal{P} \ .
\end{eqnarray}
For integer $r$ we then have 
\begin{eqnarray}
     \fl Z_n^r(\gamma,\mathcal{D})&=&
    \int \rme^{ n \sum_{\alpha=1}^r \Big\{ \rmi\sum_{\Delta = \pm 1} \int \hat{\mathcal{P}}_{\alpha}(\Delta,t,h,\gamma)\mathcal{P}_{\alpha}(\Delta,t,h,\gamma)\rmd t \rmd h - \gamma \mathcal{E}\big[\mathcal{P}_{\alpha} (.)\big] \Big\} }\times \nonumber \\
    \fl  &\times& \bigg\{\int \prod_{\alpha=1}^r  \rme^{- \rmi \sum_{i=1}^n \hat{\mathcal{P}}_{\alpha}(\Delta_i,T_i,\mathbf{X}_i'\bbeta_{\alpha}) -\gamma r(\bbeta_{\alpha}) }\rmd \bbeta_{\alpha}\bigg\}\  \prod_{\alpha=1}^r \mathcal{D}\hat{\mathcal{P}}_{\alpha} \ \mathcal{D}\mathcal{P}_{\alpha} \ .
\end{eqnarray}
Taking the expectation with respect to the data-set 
\begin{equation}
     \mathbb{E}_{\mathcal{D}}\Big[Z_n^r(\gamma,\mathcal{D})\Big]=
    \int \rme^{ n \sum_{\alpha=1}^r \mathcal{A}\big[\hat{\mathcal{P}}_{\alpha},\mathcal{P}_{\alpha}\big] }\mathcal{W}\big[\big\{\hat{\mathcal{P}}_{\alpha}\big\}\big]
      \prod_{\alpha=1}^r \mathcal{D}\hat{\mathcal{P}}_{\alpha} \ \mathcal{D}\mathcal{P}_{\alpha}
\end{equation}
where 
\begin{equation}
    \mathcal{A}\big[\gamma,\hat{\mathcal{P}}_{\alpha},\mathcal{P}_{\alpha}\big] =  \rmi\sum_{\Delta = \pm 1} \int \hat{\mathcal{P}}_{\alpha}(\Delta,t,h)\mathcal{P}_{\alpha}(\Delta,t,h)\rmd t \rmd h - \gamma \mathcal{E}\big[\mathcal{P}_{\alpha} (.)\big]
\end{equation}
and
\begin{equation}
   \fl  \mathcal{W}\big[\gamma,\big\{\hat{\mathcal{P}}_{\alpha}\big\}\big] = \int \Big(\mathbb{E}_{\Delta,T,\mathbf{X}}\Big[ \rme^{- \rmi \sum_{\alpha=1}^r\hat{\mathcal{P}}_{\alpha}(\Delta,T,\mathbf{X}'\bbeta_{\alpha}) }\Big]\Big)^n \rme^{-\gamma r(\bbeta_{\alpha})} \prod_{\alpha=1}^r \rmd \bbeta_{\alpha} \ .
\end{equation}
We notice that the expression above depends on $\bbeta_{\alpha}$ only via the linear predictors $\mathbf{X}_i' \bbeta_{\alpha}$.
Since $\mathbf{X}\sim \mathcal{N}(\bm{0},\frac{1}{p}\bm{I}_{p\times p})$, we have that
\begin{equation}
    \mathbf{Y} = (Y_{0},Y_{1},\dots,Y_r)\sim \mathcal{N}\big(\bm{0}, \mathbf{C}(\{\bbeta_{\alpha}\})\big), \quad Y_{\alpha} := \mathbf{X}' \bbeta_{\alpha}
\end{equation}
with 
\begin{equation}
   \mathbf{C}(\{\bbeta_{\alpha}\})\big) = \left( \begin{array}{cc}
        \theta_0^2 & \mathbf{M}'\\
        \mathbf{M} & \mathbf{R} 
    \end{array} \right)
\end{equation}
where 
\begin{eqnarray}
    \theta_0^2 &=& \|\bbeta_0\|^2/p \\
     \mathbf{M} &=& (M_{\alpha})_{\alpha=1}^r, \quad M_{\alpha} := \bbeta_0' \bbeta_{\alpha}/p\\
     \mathbf{R} &=& (R_{\alpha,\rho})_{\alpha,\rho=1}^r, \quad R_{\alpha,\rho} := \bbeta_{\alpha}' \bbeta_{\rho}/p  = R_{\rho,\alpha} \ .
\end{eqnarray}
Then, introducing a matrix delta function, we can write 
\begin{equation}
     \mathcal{W}\big[\big\{\hat{\mathcal{P}}_{\alpha}\big\}\big] = \int  \rme^{ip{\rm Tr}(\hat{\mathbf{C}}\mathbf{C}) + p\phi(\hat{\mathbf{C}}) + n \varphi(\gamma,\mathbf{C})}  \ \rmd \hat{\mathbf{C}}\rmd \mathbf{C}
\end{equation}
with 
\begin{eqnarray}
     \varphi[\gamma,\mathbf{C},\{\hat{\mathcal{P}}_{\alpha}(.)\}_{\alpha=1}^r] = \log \mathbb{E}_{\Delta,T,\mathbf{Y}}^{C}\Big[ \rme^{- \rmi \sum_{\alpha=1}^r\hat{\mathcal{P}}_{\alpha}(\Delta,T,Y_{\alpha}) }\Big]\\
     \phi(\gamma,\hat{\mathbf{C}}) = \frac{1}{p}\log \int \rme^{-\rmi p{\rm Tr}(\hat{\mathbf{C}}\mathbf{C}(\{\bbeta_{\alpha}\})) -\gamma \frac{1}{2} \rmr(\bbeta_{\alpha})}\prod_{\alpha=1}^r \rmd \bbeta_{\alpha} \ .
\end{eqnarray}
Putting everything together, we have obtained a so-called \say{saddle point} (functional) integral 
\begin{equation}
    \mathbb{E}_{\mathcal{D}}\Big[Z_n^r(\gamma,\mathcal{D})\Big] = \int \rme^{ -n \psi\big[\{\mathcal{P}_{\alpha},\hat{\mathcal{P}}_{\alpha}\}_{\alpha=1}^r,\hat{\mathbf{C}},\mathbf{C}\big] } \prod_{\alpha=1}^r  \ \mathcal{D}\hat{\mathcal{P}}_{\alpha} \ \mathcal{D}\mathcal{P}_{\alpha} \rmd \mathbf{C} \rmd \hat{\mathbf{C}}
\end{equation}
where
\begin{eqnarray}
\label{def : psi}
    -\psi\big[\{\mathcal{P}_{\alpha},\hat{\mathcal{P}}_{\alpha}\}_{\alpha=1}^r,\hat{\mathbf{C}},\mathbf{C}\big]&=& \sum_{\alpha=1}^r   \mathcal{A}\big[\hat{\mathcal{P}}_{\alpha},\mathcal{P}_{\alpha}\big]  + \varphi[\gamma,\mathbf{C},\{\hat{\mathcal{P}}_{\alpha}(.)\}_{\alpha=1}^r]+\nonumber \\
    &+&\rmi \zeta {\rm Tr}(\hat{\mathbf{C}}\mathbf{C}) + \zeta \phi(\gamma,\hat{\mathbf{C}})  \ . 
\end{eqnarray}

\subsection{Saddle point integration}
The idea is now to evaluate the integral via the saddle point method by interchanging the limits $n\rightarrow \infty$ and $r\rightarrow 0$, as customary in these calculations \cite{Coolen_17, Coolen_2020, Loureiro_2022, Massa_24}. 
We first derive the stationary conditions with respect to the functions $\mathcal{P}_{\alpha}$ and $\hat{\mathcal{P}}_{\alpha}$, which are obtained via functional differentiation
\begin{eqnarray}
     \mathcal{P}_{\alpha}(\Delta,t,h,\gamma)  &=& \frac{ \mathbb{E}^{\mathbf{C}}_{\mathbf{Y}}\Big[f(\Delta,t|Y_0)\delta(h - Y_{\alpha})\rme^{- \rmi \sum_{\alpha=1}^r  \hat{\mathcal{P}}_{\alpha}(\Delta,t,Y_{\alpha},\gamma)}\Big]}{ \mathbb{E}^{\mathbf{C}}_{\Delta,T,\mathbf{Y}}\Big[\rme^{- \rmi  \sum_{\alpha=1}^r  \hat{\mathcal{P}}_{\alpha}(\Delta,T,Y_{\alpha},\gamma)}\Big]} \\
     \rmi \hat{\mathcal{P}}_{\alpha}(\Delta,t,h,\gamma)  &=& \gamma  \Delta\Big[ \log \mathcal{S}_{\alpha}(t,\gamma)   - h \Big]   + \gamma \rme^{h}   \Lambda_{\alpha}(t,\gamma) \ . 
\end{eqnarray}
where
\begin{eqnarray}
     \mathcal{S}_{\alpha}(t,\gamma) &=& \sum_{\Delta' = \pm 1} \int \Theta(t-t') \rme^{h'} \mathcal{P}_{\alpha}(\Delta',t,'h')  \rmd h'\rmd t' \\
     \Lambda_{\alpha}(t,\gamma) &:=& \sum_{\Delta' = \pm 1}\int \frac{\Delta' \Theta(t-t')\mathcal{P}_{\alpha}(\Delta',t',h',\gamma) }{\mathcal{S}_{\alpha}(t',\gamma) } \rmd t' \rmd h' \ .
\end{eqnarray}
Equivalently 
\begin{eqnarray}
    \fl && \mathcal{P}_{\alpha}(\Delta,t,h,\gamma)  = \nonumber\\
    \fl && \frac{ \mathbb{E}^{\mathbf{C}}_{\mathbf{Y}}\Big[f(\Delta,t|Y_0)\delta(h - Y_{\alpha})\rme^{- \sum_{\alpha=1}^r  \gamma \big( \Delta \log \mathcal{S}_{\alpha}(t,\gamma) + \exp(Y_{\alpha})\Lambda_{\alpha}(t,\gamma) - \Delta Y_{\alpha}\big)}\Big]}{ \mathbb{E}^{\mathbf{C}}_{\Delta,T,\mathbf{Y}}\Big[\rme^{- \sum_{\alpha=1}^r  \gamma \big( \Delta \log \mathcal{S}_{\alpha}(T,\gamma) + \exp(Y_{\alpha})\Lambda_{\alpha}(T,\gamma) - \Delta Y_{\alpha}\big)}\Big]} \ .
\end{eqnarray}
Furthermore at the saddle point we have 
\begin{eqnarray}
   \fl  &&\varphi[\gamma,\mathbf{C},\{\hat{\mathcal{P}}_{\alpha}(.)\}_{\alpha=1}^r] = \nonumber \\
   \fl &&\varphi(\gamma,\mathbf{C}) = \log \mathbb{E}^{\mathbf{C}}_{\Delta,T,\mathbf{Y}}\Big[\rme^{- \sum_{\alpha=1}^r  \gamma \big( \Delta \log \mathcal{S}_{\alpha}(T, \gamma) + \exp(Y_{\alpha})\Lambda_{\alpha}(T,\gamma) - \Delta Y_{\alpha}\big)}\Big] \nonumber
\end{eqnarray}
and hence 
\begin{eqnarray}
     \fl -\tilde{\psi}(\gamma,\hat{\mathbf{C}},\mathbf{C}) &:=& - \underset{\{\mathcal{P}_{\alpha},\hat{\mathcal{P}}_{\alpha}\}}{\rm extr} \psi = 
     \rmi \zeta {\rm Tr}(\hat{\mathbf{C}}\mathbf{C}) + \zeta \phi(\gamma,\hat{\mathbf{C}})+ \varphi(\gamma,\mathbf{C}) + \nonumber \\
      \fl &+&  \gamma\sum_{\Delta = \pm 1} \int \rme^{h}  \sum_{\alpha = 1}^r  \Lambda_{\alpha}(t,\gamma)  \mathcal{P}_{\alpha}(\Delta,t,h,\gamma)  \rmd t \rmd h 
    % &+&   \gamma \sum_{\rho=1}^r \mathbb{E}^{\mathbf{C}}_{\Delta,T,\mathbf{Y}} \Big[\rme^{Y_{\rho}}\Lambda_{\rho}(t,\gamma)\rme^{- \sum_{\alpha=1}^r  \gamma \big( \Delta \log \mathcal{S}_{\alpha}(\Delta,T) + \exp(Y_{\alpha})\Lambda_{\alpha}(T) - \Delta Y_{\alpha}\big)}\Big]   \ . \nonumber 
\end{eqnarray}
With a modest amount of foresight we take the following change of variables 
\begin{equation}
    \rmi \hat{\mathbf{C}} = \frac{1}{2} \mathbf{D}
\end{equation}
which is expected from previous similar calculations and aids the book-keeping.
In principle we could now derive the saddle point equations for the elements of the matrices  $\mathbf{C}$ nor $\mathbf{D}$ and then take the limit $r\rightarrow 0$.
In practice we will assume the replica symmetric ansatz 
\begin{equation}
   \fl  \mathbf{C} = \left( \begin{array}{ccccc}
        \theta_0^2  &m   & \dots  &\dots  & m \\
        m   &\rho    & q & \dots &q \\
        \vdots   &q   & \rho & \ddots & \vdots \\
        \vdots   &\vdots    & \ddots  & \ddots  & q \\
        m &q & \dots  &q & \rho 
    \end{array}\right)\quad
    \mathbf{D} =  \left( \begin{array}{ccccc}
       0   &\hat{m}   & \dots  &\dots  & \hat{m} \\
        \hat{m}   &\hat{\rho}    & -\hat{q} & \dots & -\hat{q} \\
        \vdots   &-\hat{q}   & \hat{\rho} & \ddots & \vdots \\
        \vdots   &\vdots    & \ddots  & \ddots  & -\hat{q} \\
       \hat{m} & -\hat{q} & \dots  &-\hat{q} & \hat{\rho} 
    \end{array}\right)
\end{equation}
directly from now on.
Note that this directly implies that 
\begin{equation}
    \mathbf{C}^{-1} = \left( \begin{array}{ccccc}
        \tilde{\mu}_r   &\tilde{m}_r   & \dots  &\dots  & \tilde{m}_r \\
        \tilde{m}_r   &\tilde{\rho}_r    & \tilde{q}_r & \dots &\tilde{q}_r \\
        \vdots   &\tilde{q}_r   & \tilde{\rho}_r & \ddots & \vdots \\
        \vdots   &\vdots    & \ddots  & \ddots  & \tilde{q}_r \\
        \tilde{m}_r &\tilde{q}_r & \dots  &\tilde{q}v & \tilde{\rho}_r 
    \end{array}\right) \ .
\end{equation}
After some algebra one obtains that
\begin{eqnarray}
    \tilde{\mu}_r &=& \frac{\rho+q(r-1)}{\theta_0^2(\rho +q(r-1))-rm^2} \\
    \tilde{m}_r &=& \frac{m}{rm^2-\theta_0^2(\rho + q(r-1))}\\ 
    \tilde{q}_r &=& \frac{1}{\rho-q} \; \frac{m^2-q\theta_0^2}{\theta_0^2(\rho + q(r-1))-rm^2} \\
    \tilde{\rho}_r &=& \frac{1}{\rho-q}\Big(1+\frac{m^2-q\theta_0^2}{\theta_0^2(\rho + q (r-1))-rm^2}\Big) \ .
\end{eqnarray}
We are now in position to compute the \say{limit} $r\rightarrow 0$, i.e. to obtain the replica symmetric symmetric free energy.

\subsubsection{Simplification of $\phi$ within RS ansatz}
Let's remind the definition of the potential $\phi$
\begin{equation}
    \phi(\gamma,\mathbf{D}) = \frac{1}{p}\log \int \rme^{-\frac{1}{2}p{\rm Tr}(\mathbf{D}\mathbf{C}(\{\bbeta_{\alpha}\})) -\gamma \rmr(\bbeta_{\alpha})}\prod_{\alpha=1}^r \rmd \bbeta_{\alpha} \ .
\end{equation}
Using the replica symmetric ansatz we get 
\begin{eqnarray}
    \fl  \frac{1}{2} p{\rm Tr}(\mathbf{D}\mathbf{C}(\{\bbeta_{\alpha}\})) &=& \sum_{\rho=1}^r \big( \hat{m}\bbeta_{0}'\bbeta_{\rho}  + \frac{1}{2} \hat{\rho} \bbeta_{\rho}'\bbeta_{\rho}\big) - \sum_{\rho,\alpha\neq \rho}^r \hat{q} \bbeta_{\rho}\bbeta_{\alpha} \nonumber \\
    \fl &=&    \sum_{\rho=1}^r \big( \hat{m}\bbeta_{0}'\bbeta_{\rho}  + \frac{1}{2} (\hat{\rho} +\hat{q})  \|\bbeta_{\rho}\|^2\big) - \frac{1}{2}\hat{q}\big\|\sum_{\rho=1}^r  \bbeta_{\rho}\big\|^2 \nonumber
\end{eqnarray}
and via Gaussian linearization
\begin{equation}
     \fl \phi^{(r)}_{RS}(\gamma,\hat{m},\hat{q},\hat{\rho}) =  \frac{1}{p}\log \mathbb{E}_{\mathbf{Z}}\bigg[\bigg(\int \rme^{-\frac{1}{2}(\hat{\rho} +\hat{q}) \|\mathbf{x}\|^2 - (\hat{m}\bbeta_{0} + \sqrt{\hat{q}} \mathbf{Z} )' \mathbf{x} -\gamma \rmr(\bx)}\rmd \mathbf{x}\bigg)^r\bigg] \nonumber \ .
\end{equation}
Since we are finally interested in taking the limit $r\rightarrow 0$ it is convenient to expand the integrand for small $r$, thus obtaining 
\begin{eqnarray}
    \fl\phi^{(r)}_{RS}(\gamma,\hat{m},\hat{q},\hat{\rho}) =  r\frac{1}{p}\mathbb{E}_{\mathbf{Z}}\bigg[\log \int \rme^{-\frac{1}{2}(\hat{\rho} +\hat{q}) \|\mathbf{x}\|^2 - (\hat{m}\bbeta_{0} + \sqrt{\hat{q}} \mathbf{Z} )' \mathbf{x} -\gamma \rmr(\bx)}\rmd \mathbf{x}\bigg]  + O(r^2) \ . \nonumber
\end{eqnarray}
Since the integrand in the expression above factorizes over the components of $\bx$ and $\mathbf{Z}$, this can be further re-written in terms of the distribution of the entries of $\bbeta_0$ 
\begin{equation}
\fl\phi^{(r)}_{RS}(\gamma,\hat{m},\hat{q},\hat{\rho}) =  r\mathbb{E}_{\beta_0, Z}\bigg[\log \int \rme^{-\frac{1}{2}(\hat{\rho} +\hat{q}) x^2 - (\hat{m}\beta_0 + \sqrt{\hat{q}} Z )x -\gamma \rmr(x)}\rmd x\bigg]  + O(r^2) \ .\nonumber
\end{equation}
\subsubsection{Simplification of $\varphi$ within RS ansatz}
Inserting the replica symmetric ansatz, we obtain
\begin{equation}
   \fl f(y_0,y_1,\dots,y_{r}) \propto \exp\Big\{-\Big(\frac{1}{2}\tilde{\mu}y_0^2 +  \sum_{\rho=1}^r \big(\tilde{m} y_0 y_{\rho} +\frac{1}{2} (\tilde{\rho} - \tilde{q})y_{\rho}^2\big) + \frac{1}{2} \tilde{q}\big(\sum_{\rho=1}^r y_{\rho}\big)^2 \Big)\Big\}
\end{equation}
and via Gaussian linearization 
\begin{equation}
\label{measure_simplification}
   \fl f(y_0,y_1,\dots,y_{r}) \propto \mathbb{E}_{Q}\Big[\exp\Big\{ -\frac{1}{2}\tilde{\mu}y_0^2 -  \sum_{\rho=1}^r\big[ \frac{1}{2} (\tilde{\rho} - \tilde{q})y_{\rho}^2  + \big(\tilde{m} y_0 + \rmi \sqrt{\tilde{q}}Q \big)y_{\rho}\big] \Big\}\Big]
\end{equation}
with $Q\sim \mathcal{N}(0,1)$.
Upon setting $\sqrt{\tilde{\mu}}y_0 = z_0$ we obtain 
\begin{eqnarray}
     \fl &&\varphi^{(r)}_{RS}(\gamma,\tilde{\mu}_r,\tilde{m}_r,\tilde{q}_r,\tilde{\rho}_r) \nonumber \\
     \fl && = \log  \frac{\mathbb{E}_{\Delta,T,Z_0,Q}\Bigg[\bigg(\int \rme^{ -\frac{1}{2} (\tilde{\rho} - \tilde{q})x^2  - \big(\tilde{m}/\sqrt{\tilde{\mu}} Z_0 + \rmi \sqrt{\tilde{q}}Q \big)x - \gamma g(x,\Lambda^{(r)}(T,\gamma),\Delta)}\frac{\rmd x}{\sqrt{2\pi}}\bigg)^r\Bigg]}{\mathbb{E}_{\Delta,T,Z_0,Q}\Bigg[\bigg(\int \rme^{ -\frac{1}{2} (\tilde{\rho} - \tilde{q})x^2  - \big(\tilde{m}/\sqrt{\tilde{\mu}} Z_0 + \rmi \sqrt{\tilde{q}}Q \big)x}\bigg)^r\Bigg]} \ . \nonumber 
\end{eqnarray}
with $Z_0 \sim \mathcal{N}(0,1)$, $Z_0\perp Q$ and where we have defined 
\begin{equation}
    g(x,y,z) =  \exp(x)y - z x  \ .
\end{equation}
\subsection{Replica symmetric functional equation}
Using (\ref{measure_simplification}) we obtain 
\begin{eqnarray}
\label{def : rs_P_gamma_r}
     \fl &&\mathcal{P}^{(r)}(\Delta,t,h,\gamma)  =  \bigg(\mathbb{E}_{\Delta,T,Z_0,Q}\Big[ \mathcal{Z}^{r}(\gamma,\Delta,T,Z_0,Q)\Big]\bigg)^{-1}\times   \nonumber \\
     \fl && \times \mathbb{E}_{Z_0,Q}\Bigg[\frac{f(\Delta,t|Z_0/\sqrt{\tilde{\mu}})\int \delta(h-x) \ \rme^{-\Big\{\frac{1}{2} (\tilde{\rho} - \tilde{q})x^2  + \big(\tilde{m}/\sqrt{\tilde{\mu}}Z_0 + \rmi \sqrt{\tilde{q}}Q \big)x + \gamma g(x,\Lambda^{(r)}(T,\gamma),\Delta)\Big\}} \rmd x }{\bigg(\int \rme^{-\Big\{\frac{1}{2} (\tilde{\rho} - \tilde{q})x^2  + \big(\tilde{m}/\sqrt{\tilde{\mu}}Z_0 + \rmi \sqrt{\tilde{q}}Q \big)x + \gamma g(x,\Lambda^{(r)}(T,\gamma),\Delta)\Big\}}\rmd x\bigg)^{1-r}}\Bigg]\nonumber 
\end{eqnarray}
where we used the shorthand 
\begin{equation}
   \fl \mathcal{Z}(\gamma,\Delta,T,Z_0,Q) = \int \rme^{-\gamma\Big\{\frac{1}{2} \frac{(\tilde{\rho} - \tilde{q})}{\gamma}x^2  + \frac{1}{\gamma}\big(\tilde{m}/\sqrt{\tilde{\mu}}Z_0 + \rmi \sqrt{\tilde{q}}Q \big)x +\gamma g(x,\Lambda^{(r)}(T,\gamma),\Delta)\Big\}}\rmd x \ .
\end{equation}

Let us stop and recall what we have achieved so far. By assuming the RS ansatz we have obtained 
\begin{equation}
    -\lim_{n \rightarrow\infty }\frac{1}{n} \log \mathbb{E}_{\mathcal{D}}\Big[Z_n^r(\gamma,\mathcal{D})\Big] = \underset{m,q,\rho,\hat{\mu},\hat{m},\hat{q},\hat{\rho}}{{\rm extr}}\tilde{\psi}^{(r)}_{RS}(m,q,\rho,\hat{m},\hat{q},\hat{\rho})
\end{equation}
with
\begin{eqnarray}
     \fl &&-\tilde{\psi}^{(r)}_{RS}(\dots ) =  \\
     \fl && \zeta (\mu\hat{\mu} + r m \hat{m} + r (\rho \hat{\rho}- q \hat{q}) + r^2 q\hat{q}\big) + \zeta \phi^{(r)}_{RS}(\gamma,\hat{m},\hat{q},\hat{\rho}) + \varphi^{(r)}_{RS}(\gamma,\tilde{\mu}_r,\tilde{m}_r,\tilde{q}_r,\tilde{\rho}_r)   + \nonumber\\
    \fl && \gamma r\mathbb{E}_{\Delta,T,Z_0,Q} \Bigg[\Lambda^{(r)}(t,\gamma)\frac{\int \rme^{x} \rme^{-\Big\{\frac{1}{2} (\tilde{\rho} - \tilde{q})x^2  + \big(\tilde{m}/\sqrt{\tilde{\mu}}Z_0 + \rmi \sqrt{\tilde{q}}Q \big)x + \gamma g(x,\Lambda^{(r)}(T,\gamma),\Delta)\Big\}}\rmd x}{\bigg(\int  \rme^{-\Big\{\frac{1}{2} (\tilde{\rho} - \tilde{q})x^2  + \big(\tilde{m}/\sqrt{\tilde{\mu}}Z_0 + \rmi \sqrt{\tilde{q}}Q \big)x + \gamma g(x,\Lambda^{(r)}(T,\gamma),\Delta)\Big\}}\rmd x\bigg)^{1-r}}\Bigg]   \ .\nonumber
\end{eqnarray}
where now
\begin{eqnarray}
     \fl &&\mathcal{S}^{(r)}(t',\gamma) = \\
     \fl &&=\mathbb{E}_{\Delta,T,Z_0,Q} \Bigg[\Theta(T-t')\frac{\int \rme^{x} \rme^{-\Big\{\frac{1}{2} (\tilde{\rho} - \tilde{q})x^2  + \big(\tilde{m}/\sqrt{\tilde{\mu}}Z_0 + \rmi \sqrt{\tilde{q}}Q \big)x + \gamma g(x,\Lambda^{(r)}(T,\gamma),\Delta)\Big\}}\rmd x}{\bigg(\int  \rme^{-\Big\{\frac{1}{2} (\tilde{\rho} - \tilde{q})x^2  + \big(\tilde{m}/\sqrt{\tilde{\mu}}Z_0 + \rmi \sqrt{\tilde{q}}Q \big)x + \gamma g(x,\Lambda^{(r)}(T,\gamma),\Delta)\Big\}}\rmd x\bigg)^{1-r}}\Bigg] \nonumber\\
    \fl && \Lambda^{(r)}(t,\gamma) := \sum_{\Delta' = \pm 1 }\int \frac{\Delta' \Theta(t-t') }{\mathcal{S}(t',\gamma)} \mathcal{P}^{(r)}(\Delta', t', h', \gamma) \rmd t' \rmd h' \nonumber \ . 
\end{eqnarray}

\subsection{The limit $r\rightarrow 0 $}
Taking the limit $r\rightarrow 0$, we get the replica symmetric free energy 
\begin{equation}
     f_{RS}(\gamma) = \underset{\mu,m,q,\rho,\hat{m},\hat{q},\hat{\rho}}{{\rm extr}}\tilde{f}_{RS}(\gamma,\mu,m,q,\rho,\hat{m},\hat{q},\hat{\rho}) 
\end{equation}
with 
\begin{equation}
    \tilde{f}_{RS}(\gamma,m,q,\rho,\hat{m},\hat{q},\hat{\rho}) := \lim_{r\rightarrow 0}\frac{1}{\gamma r} \tilde{\psi}^{(r)}_{RS}(\gamma,m,q,\rho,\hat{m},\hat{q},\hat{\rho}) \ .
\end{equation}
where
\begin{eqnarray}
     \fl &&-\tilde{f}_{RS}(\dots) = \zeta (  m \hat{m} + r(\rho \hat{\rho}- q \hat{q}) \big) +  \zeta \tilde{\phi}_{RS}(\gamma,\hat{m},\hat{q},\hat{\rho}) +  \tilde{\varphi}_{RS}(\gamma,m,q,\rho) +\nonumber \\
     \fl &&  \mathbb{E}_{\Delta,T,Z_0,Q} \Bigg[\Lambda(t,\gamma)\frac{\int \rme^{x} \rme^{ - \frac{1}{\rho-q}\Big\{\frac{1}{2}x^2  - \big((m/\theta_0)Z_0 + \sqrt{q - (m/\theta_0)^2}Q \big)x\Big\} - \gamma  g(x,\Lambda(t,\gamma),\Delta)} \rmd x}{\int \rme^{ - \frac{1}{\rho-q}\Big\{\frac{1}{2}x^2  - \big((m/\theta_0)Z_0 + \sqrt{q - (m/\theta_0)^2}Q \big)x\Big\} - \gamma  g(x,\Lambda(t,\gamma),\Delta)}\rmd x}\Bigg]\ . 
\end{eqnarray}
Above, we introduced the following definitions 
\begin{eqnarray}
      \fl &&\tilde{\phi}_{RS}(\gamma,\hat{m},\hat{q},\hat{\rho}) = \lim_{r\rightarrow 0}\frac{1}{\gamma r} \phi_{RS}^{(r)}(\gamma,\hat{m},\hat{q},\hat{\rho}) \nonumber \\
     \fl && \hspace{2 cm} = \mathbb{E}_{\beta_0, Z}\bigg[\log \int \rme^{-\frac{1}{2}(\hat{\rho} +\hat{q}) x^2 - (\hat{m}\beta_0 + \sqrt{\hat{q}} Z )x -\gamma \rmr(x)}\rmd x\bigg]\nonumber\\
    \fl  && \tilde{\varphi}_{RS}(\gamma,m,q,\rho) = \lim_{r\rightarrow 0}\frac{1}{\gamma r} \varphi^{(r)}_{RS}(\gamma,\tilde{\mu}_r,\tilde{m}_r,\tilde{q}_r,\tilde{\rho}_r)  \nonumber\\
    \fl  &&= \mathbb{E}_{\Delta,T,Z_0,Q}\Bigg[\log\frac{ \int \rme^{ - \frac{1}{\rho-q}\Big\{\frac{1}{2}x^2  - \big((m/\theta_0)Z_0 + \sqrt{q - (m/\theta_0)^2}Q \big)x\Big\} - \gamma  g(x,\Lambda(t,\gamma),\Delta)}\rmd x}{ \int \rme^{ - \frac{1}{\rho-q}\Big\{\frac{1}{2}x^2  - \big((m/\theta_0)Z_0 + \sqrt{q - (m/\theta_0)^2}Q \big)x\Big\}} \rmd x}\Bigg] \nonumber 
    %&=& \frac{1}{p}\sum_{\mu=1}^p  \mathbb{E}_{\varphi}\bigg[ \frac{1}{2 \gamma}\Big(\frac{\hat{m}\beta_{0,mu} + i \sqrt{\hat{q}}\varphi}{\sqrt{\hat{\rho} -\hat{q}}}\Big)^2\bigg] = \frac{1}{2 \gamma } \frac{\hat{m}^2 \theta_0^2 - \hat{q}}{\hat{\rho} -\hat{q}}, \quad \theta_0 := \frac{\|\bbeta_0\|^2}{p}\nonumber \\
\end{eqnarray}
and we have used that 
\begin{eqnarray}
    \tilde{\mu} &=& \lim_{r\rightarrow 0 } \tilde{\mu}_r = 1/\theta_0^2\\
    \tilde{m} &=& \lim_{r\rightarrow 0 } \tilde{m}_r =-\frac{m}{\theta_0^2(\rho-q )}\\ 
    \tilde{q} &=& \lim_{r\rightarrow 0 } \tilde{q}_r = \frac{1}{\rho-q}  \frac{m^2-q\theta_0^2}{\theta_0^2(\rho -q)} \\
    \tilde{\rho} - \tilde{q} &=& \lim_{r\rightarrow 0 }\tilde{\rho}_r - \tilde{q}_r= \frac{1}{\rho-q} \ .
\end{eqnarray}
Assuming the following scaling
\begin{equation}
    \tau =  \gamma/(\tilde{\rho} - \tilde{q}) = \gamma(\rho - q)
\end{equation}
we obtain 
\begin{eqnarray}
     \fl  &&\tilde{\varphi}_{RS}(\gamma,m,q,\rho) =  \nonumber \\
     \fl  && = \mathbb{E}_{\Delta,T,Z_0,Q}\Bigg[ \frac{1}{\gamma}\log \frac{\int \rme^{-\gamma \Big\{\frac{1}{2}\big(x -( (m/\theta_0)Z_0 + \sqrt{q - (m/\theta_0)^2}Q)\big)^2/\tau +g\big(x,\Lambda(T,\gamma),\Delta\big)\Big\}}\rmd x}{\int \rme^{-\gamma \Big\{\frac{1}{2}\big(x -( (m/\theta_0)Z_0 + \sqrt{q - (m/\theta_0)^2}Q)\big)^2/\tau\Big\}}\rmd x} \Bigg] \ .
\end{eqnarray}
Similarly, taking the additional rescaling
\begin{equation}
    1/\hat{\tau} = (\hat{\rho} +\hat{q}) /\gamma, \quad \hat{m} = \gamma \hat{m}, \quad \hat{q} = \gamma^2 \hat{q}
\end{equation}
we have 
\begin{eqnarray}
    \fl  && \tilde{\phi}_{RS}(\gamma,\hat{m},\hat{q},\hat{\rho}) =  \mathbb{E}_{\beta_0,Z}\bigg[ \frac{1}{\gamma}\log  \int \rme^{-\gamma\Big\{\frac{1}{2\hat{\tau} } x^2 + (\hat{m}\beta_0+ \sqrt{\hat{q}} Z )x +  \rmr(x)\Big\}}\rmd x\bigg] =  \nonumber \\ 
    \fl  &&=\frac{1}{2}\hat{\tau} \Big( \hat{m}^2 \theta_0^2  + \hat{q}\Big) +  \mathbb{E}_{\beta_0,Z}\bigg[  \frac{1}{\gamma}\log \int \rme^{-\gamma \Big\{\frac{1}{2}\frac{1}{\hat{\tau}} \big\{\tilde{x} + \hat{\tau} (\hat{m}\beta_0 +  \sqrt{\hat{q}} Z)\big\}^2  + \rmr(\tilde{x})\Big\} }\rmd \tilde{x}\bigg] +{\rm const}\nonumber 
\end{eqnarray}
Furthermore, in the limit $r\rightarrow 0$, the functional equation reads
\begin{eqnarray}
    \label{funct_eq_gamma}
     \fl &&\mathcal{P}(\Delta,t,h,\gamma) = \\
     \fl &&  \mathbb{E}_{Z_0,Q}\Bigg[\frac{f(\Delta,t|Z_0/\sqrt{\tilde{\mu}})\int \delta(h-x) \ \rme^{-\gamma \Big\{\frac{1}{2}\big(x -( (m/\theta_0)Z_0 + \sqrt{q - (m/\theta_0)^2}Q)\big)^2/\tau +g\big(x,\Lambda(T,\gamma),\Delta\big)\Big\}} \rmd x }{\int \rme^{-\gamma \Big\{\frac{1}{2}\big(x -( (m/\theta_0)Z_0 + \sqrt{q - (m/\theta_0)^2}Q)\big)^2/\tau +g\big(x,\Lambda(T,\gamma),\Delta\big)\Big\}}\rmd x }\Bigg]\nonumber 
\end{eqnarray}
where 
\begin{eqnarray}
    \fl  &&\mathcal{S}(t',\gamma) = \nonumber \\
    \fl && \mathbb{E}_{\Delta,T,Z_0,Q} \Bigg[\Theta(T-t')\frac{\int \rme^{x} \rme^{-\gamma \Big\{\frac{1}{2}\big(x -( (m/\theta_0)Z_0 + \sqrt{q - (m/\theta_0)^2}Q)\big)^2/\tau +g\big(x,\Lambda(T,\gamma),\Delta\big)\Big\}}\rmd x}{\int  \rme^{-\gamma \Big\{\frac{1}{2}\big(x -( (m/\theta_0)Z_0 + \sqrt{q - (m/\theta_0)^2}Q)\big)^2/\tau +g\big(x,\Lambda(T,\gamma),\Delta\big)\Big\}}\rmd x}\Bigg] \nonumber\\
    \fl &&\Lambda(t,\gamma) := \mathbb{E}_{\Delta',T'}\Bigg[\frac{\Delta' \Theta(t-T') }{\mathcal{S}(T',\gamma)}\Bigg] \nonumber \ . 
\end{eqnarray}
In the last expression we have used that 
\begin{equation}
    \int \mathcal{P}(\Delta,t,h,\gamma) \rmd h = f(\Delta, t|Z_0/\sqrt{\tilde{\mu}}) \ .
\end{equation}

\subsection{The limit $\gamma \rightarrow \infty $}
Via Laplace integration we see that 
\begin{equation}
\label{laplace_moureau}
    - \lim_{ \gamma \rightarrow \infty} \frac{1}{\gamma } \log \int \rme^{-\gamma \Big\{\frac{1}{2 \alpha} (z - x)^2  + b(z)\Big\}} \rmd z =  \mathcal{M}_{b(.)}(x,\alpha)
\end{equation}
where $\mathcal{M}_{b(.)}$ is the Moureau envelope of a convex function $b:\mathbb{R} \rightarrow \mathbb{R}$, which is defined as 
\begin{equation}
    \mathcal{M}_{b(.)}(x,\alpha) = \underset{z}{\min}\Big\{\frac{1}{2\alpha} (z-x)^2  + b(z)\Big\} \ . 
\end{equation}
Using the \say{Laplace-Moureau} identity above (\ref{laplace_moureau}), we obtain 
\begin{eqnarray}
    \fl  &&\lim_{\gamma \rightarrow \infty } \tilde{\phi}_{RS}(\gamma,\hat{m},\hat{q},\hat{\rho}) = \frac{1}{2}\hat{\tau} \Big( \hat{m}^2 \theta_0^2  + \hat{q}\Big) - \mathbb{E}_{Z, \beta_0}\Big[\mathcal{M}_{\rmr(.)}\Big(\hat{m}\beta_0 +  \sqrt{\hat{q}} Z, \hat{\tau}\Big)\Big] \nonumber \\
    \fl  && \lim_{\gamma \rightarrow \infty } \tilde{\varphi}_{RS}(\gamma,m,q,\rho) =  -\mathbb{E}_{\Delta,T,Z_0,Q}\Big[\mathcal{M}_{\tau g(.,\Lambda(T),\Delta)}\Big(m/\theta_0Z_0 + \sqrt{q - (m/\theta_0)^2}Q, \tau\Big)\Big] \nonumber  \ .
\end{eqnarray}
In the limit $\gamma\rightarrow \infty$, the functional equation (\ref{funct_eq_gamma}) reads
\begin{eqnarray}
    \label{def : P_rs}
     \mathcal{P}(\Delta,t,h)  &=& \mathbb{E}_{Z_0,Z}\Big[f(\Delta,t|Z_0)\delta\Big(h - \xi(T, Z_0,Z)\Big)\Big] \\
    \label{def: prox_cox}
     \xi(T, Z_0,Z) &:=& {\rm prox}_{g(.,\Lambda(T),\Delta)}\Big((m/\theta_0)Z_0 + \sqrt{q - (m/\theta_0)^2}Z, \tau \Big), 
\end{eqnarray}
where we used the Laplace integration method to conclude that for any \say{well behaved} function $a:\mathbb{R}^d\rightarrow\mathbb{R} $
\begin{equation}
     a\Big({\rm prox}_{b(.)}(x, \alpha)\Big) = \lim_{ \gamma \rightarrow \infty}  \log \frac{\int \rme^{-\gamma \Big\{\frac{1}{2\alpha} (z - x)^2  + b(z)\Big\}} a(z)\rmd z}{\int \rme^{-\gamma \Big\{\frac{1}{2\alpha} (z - x)^2  + b(z)\Big\}} \rmd z} 
\end{equation}
with ${\rm prox}_{b(.)}$ the proximal mapping of the convex function $b$ defined as 
\begin{equation}
   {\rm prox}_{b(.)}(x, \alpha) =  \underset{z}{\arg \min}\Big\{\frac{1}{2\alpha} (z - x)^2  + b(z)\Big\} \ .
\end{equation}
We also have 
\begin{eqnarray}
     \mathcal{S}(t) &=&\mathbb{E}_{\Delta,T,Z_0,Z}\Big[ \Theta(t-T)\rme^{\xi}\Big] \nonumber\\
     \Lambda(t) &:=& \mathbb{E}_{\Delta',T'}\Bigg[\frac{\Delta' \Theta(t-T') }{\mathcal{S}(T')}\Bigg] \nonumber 
\end{eqnarray}
where $\xi$ is defined in (\ref{def: prox_cox}).
Furthermore 
\begin{eqnarray}
    \fl  &&\lim_{\gamma \rightarrow \infty}\mathbb{E}_{\Delta,T,Z_0,Z} \Bigg[\Lambda(t,\gamma)\frac{\int \rme^{x} \rme^{ - \frac{1}{\rho-q}\Big\{\frac{1}{2}x^2  - \big((m/\theta_0)Z_0 + \sqrt{q - (m/\theta_0)^2}Z \big)x\Big\} - \gamma  g(x,\Lambda(t,\gamma),\Delta)} \rmd x}{\int \rme^{ - \frac{1}{\rho-q}\Big\{\frac{1}{2}x^2  - \big((m/\theta_0)Z_0 + \sqrt{q - (m/\theta_0)^2}Z \big)x\Big\} - \gamma  g(x,\Lambda(t,\gamma),\Delta)}\rmd x}\Bigg] \nonumber \\
    \fl  &&= \mathbb{E}_{\Delta,T}\Bigg[\frac{\Delta \mathbb{E}_{\Delta'',T'',Y''_0,Z''}\Big[ \Theta(T''-T) \rme^{\xi''}\Big]}{\mathbb{E}_{\Delta',T',Y'_0,Z'}\Big[ \Theta(T'-T)\rme^{\xi'}\Big]}\Bigg] = \mathbb{E}\Big[\Delta\Big] = {\rm const} \nonumber .
\end{eqnarray}
Putting all together (for the last time), we have obtained that 
\begin{equation}
    \fl \lim_{\gamma\rightarrow\infty}    \tilde{f}_{RS}(\gamma,m,q,\rho,\hat{m},\hat{q},\hat{\rho})  = \mathcal{F}(m,q,\tau,\hat{m},\hat{q},\hat{\tau}) = \lim_{n\rightarrow\infty}\mathbb{E}_{\mathcal{D}}\Big[\frac{1}{n}\mathcal{PL}_n(\hat{\bbeta}_n)\Big]
\end{equation}
with 
\begin{eqnarray}
   \fl  &&\mathcal{F}(m,q,\tau,\hat{m},\hat{q},\hat{\tau})=-\zeta ( m \hat{m} + \frac{1}{2} q/\hat{\tau} -\frac{1}{2}\tau\hat{q}  \big) +\frac{1}{2} \hat{\tau}\zeta \Big( \hat{m}^2 \theta_0^2  + \hat{q}\Big) +  \nonumber \\
   \fl  &&\hspace{2cm}\zeta\mathbb{E}_{Z, \beta_0}\Big[\mathcal{M}_{\rmr(.)}\Big(\hat{m}\beta_0 +  \sqrt{\hat{q}} Z, \hat{\tau}\Big)\Big] +\\
     \fl && \hspace{2cm} \mathbb{E}_{\Delta,T,Z_0,Q}\Big[\mathcal{M}_{ g(.,\Lambda(T),\Delta)}\Big(m/\theta_0Z_0 + \sqrt{q - (m/\theta_0)^2}Q, \tau\Big)\Big] + {\rm const}
\end{eqnarray}
It is now convenient to take additional change of variables in order to obtain \say{neater} formulae. 
Let us define
\begin{equation}
    w = m/\theta_0, \quad v = q - m^2/\theta_0^2, \quad \hat{w} = -\hat{\tau} \hat{m} \theta_0, \quad \hat{v} =  \hat{\tau}\sqrt{\hat{q}}
\end{equation}
then (abusing the notation yet again)
\begin{eqnarray}
   \fl  &&\mathcal{F}(w,v,\tau,\hat{w},\hat{v},\hat{\tau})=-\frac{\zeta}{2\hat{\tau}}  \Big((w - \hat{w})^2 + v^2 + \hat{v}^2(1-\tau/\hat{\tau})\Big)  + \nonumber \\
    \fl &&  \zeta \mathbb{E}_{Z, \beta_0}\Big[\mathcal{M}_{\rmr(.)}\Big(\hat{w}\frac{\beta_0}{\theta_0} +  \hat{v}Z, \hat{\tau}\Big)\Big] + \nonumber\\
    \fl  && \hspace{2cm}   \mathbb{E}_{\Delta,T,Z_0,Q}\Big[\mathcal{M}_{ g(.,\Lambda(T),\Delta)}\Big( w Z_0 + v Q , \tau\Big)\Big] + {\rm const}
\end{eqnarray}

\section{Replica symmetric equations }
\label{app_sec:rs_eqs}
It is convenient to define
\begin{eqnarray}
    &&\varphi(\beta_0, Z) := {\rm prox}_{\rmr(.)}\Big(\hat{w}\frac{\beta_0}{\theta_0} +  \hat{v}Z, \hat{\tau}\Big)\\
    &&\xi:= {\rm prox}_{g(.,\Lambda(T),\Delta)}(wZ_0 + vQ, \tau)
\end{eqnarray}
then 
\begin{eqnarray}
     \fl \frac{\partial}{\partial \hat{w}} \mathcal{F} &=& \frac{\zeta}{\hat{\tau}}   w  \frac{\zeta}{\hat{\tau}}   \mathbb{E}_{Z, \beta_0}\Big[\beta_0 \varphi\Big]/\theta_0 \ \rightarrow\ w = \mathbb{E}_{Z, \beta_0}\Big[\beta_0 \varphi\Big]/\theta_0\\
     \fl \frac{\partial}{\partial \hat{v}} \mathcal{F} &=&  \frac{\zeta}{\hat{\tau}} \hat{v} \frac{\tau}{\hat{\tau}}
      \frac{\zeta}{\hat{\tau}}  \mathbb{E}_{Z, \beta_0}\Big[Z \varphi\Big] \ \rightarrow \ \hat{v} \frac{\tau}{\hat{\tau}} =\mathbb{E}_{Z, \beta_0}\Big[Z \varphi\Big]  \\
     \fl  \frac{\partial}{\partial \hat{\tau}} \mathcal{F} &=&  \frac{\zeta}{2\hat{\tau}^2} \Big((w -\hat{w})^2 + v^2 + \hat{v^2}( 1- 2\frac{\tau}{\hat{\tau}})\Big)-  \frac{\zeta}{2\hat{\tau}^2}   \mathbb{E}_{Z, \beta_0}\Big[\big(\varphi - \hat{w}\beta_0 -  \hat{v}Z \big)^2 \Big] \\ 
    \fl  &=& \frac{\zeta}{2\hat{\tau}^2} \Big(w^2  + v^2\Big)-  \frac{\zeta}{2\hat{\tau}^2}  \mathbb{E}_{Z, \beta_0}\Big[\varphi^2\Big] \\
     \fl \frac{\partial}{\partial w} \mathcal{F} &=& -\zeta(  w -\hat{w})/\hat{\tau}+ \frac{1}{\tau} \Big(w - \mathbb{E}_{\Delta,T,Z_0,Q}\Big[Z_0\xi\Big]\Big) \\
     \fl \frac{\partial}{\partial v} \mathcal{F} &=& -\zeta v/\hat{\tau}  +  \frac{1}{\tau} \Big(v - \mathbb{E}_{\Delta,T,Z_0,Q}\Big[Q\xi\Big]\Big) \\
     \fl \frac{\partial}{\partial \tau} \mathcal{F} &=& \zeta \frac{\hat{v}^2}{\hat{\tau}^2} - \frac{1}{\tau^2}\mathbb{E}_{\Delta,T,Z_0,Q}\Big[(\xi - wZ_0 - vQ)^2\Big] \ .
\end{eqnarray}
Hence
\begin{eqnarray}
    w  &=& \mathbb{E}_{Z, \beta_0}\Big[\beta_0 \varphi\Big]/\theta_0\\
    \hat{v}\frac{\tau}{\hat{\tau}} &=&\mathbb{E}_{Z, \beta_0}\Big[Z \varphi\Big] \\
    (w^2+v^2) &=&  \mathbb{E}_{Z, \beta_0}\Big[\varphi^2\Big]\\
     \hat{w} &=& w - \frac{\hat{\tau}}{\zeta \tau}\Big(w -\mathbb{E}_{\Delta,T,Z_0,Q}\Big[Z_0\xi\Big]\Big)\\
    v(1 -\zeta\tau/\hat{\tau}) &=&\mathbb{E}_{\Delta,T,Z_0,Q}\Big[Q\xi\Big] \\
    \zeta \hat{v}^2 &=& \frac{\hat{\tau}^2}{\tau^2}\mathbb{E}_{\Delta,T,Z_0,Q}\Big[(\xi - wZ_0 - vQ)^2\Big]
\end{eqnarray}

\section{Explicit computation for the lasso regularization}
\label{app_sec : computation_lasso}
Here we assume that 
\begin{equation}
    \mathbf{e}_{\mu}'\bbeta_0 \sim \mathcal{N}(0, \theta_0^2/\nu ), \quad \nu := s/p \ .
\end{equation}
For the Lasso regularization
\begin{eqnarray}
    \prox_{\rmr(.)}(x, \hat{\tau}) &=&  {\rm st}(x, \alpha\hat{\tau}), \\
     {\rm st}(x, \alpha) &:=& {\rm relu}(x-\alpha) - {\rm relu}(-x-\alpha)%\bm{1}[x>\alpha] (x-\alpha) + \bm{1}[x<-\alpha] (x+\alpha)
\end{eqnarray}
with 
\begin{equation}
    {\rm relu}(x) = x\Theta(x) \ .
\end{equation}
Before proceeding, let us remember some useful fact 
\begin{eqnarray}
    \mathbb{E}_Z\Big[\Theta(\sigma Z + \mu -\alpha)\Big] &=& \Phi\big(\frac{\alpha-\mu}{\sigma}\big) \\%\frac{1}{2} - \frac{1}{2} {\rm erf}\big(\frac{\alpha}{\sqrt{2}}\big)\\
    \mathbb{E}_Z\Big[Z\Theta(\sigma Z +\mu -\alpha)\Big] &=& \frac{\rme^{-\frac{1}{2}\big(\frac{\alpha-\mu}{\sigma}\big) ^2}}{\sqrt{2\pi}}\\
     \mathbb{E}_Z\Big[Z^2\Theta(\sigma Z + \mu -\alpha)\Big] &=&   \Phi\big(\frac{\alpha-\mu}{\sigma}\big) +  \frac{\alpha-\mu}{\sigma} \frac{\rme^{-\frac{1}{2}\big(\frac{\alpha-\mu}{\sigma}\big)^2 }}{\sqrt{2\pi}} %\frac{1}{2} - \frac{1}{2} {\rm erf}\big(\frac{\alpha}{\sqrt{2}}\big) + \frac{\alpha\rme^{-\frac{1}{2}\alpha^2}}{\sqrt{2\pi}} \ .
\end{eqnarray}
where 
\begin{equation}
    \Phi(x) := \int_{x}^{+\infty} \frac{1}{\sqrt{2\pi}}\rme^{-\frac{1}{2} x^2} \rmd x %= \frac{1}{2} - \frac{1}{2} {\rm erf}(x)
\end{equation}
The identities above imply
\begin{equation}
    \mathbb{E}_Z[Z {\rm relu}(\sigma Z + \mu - \alpha)] = \sigma \Phi\big(\frac{\alpha-\mu}{\sigma}\big) 
\end{equation}
and hence
\begin{equation}
    \mathbb{E}_Z[Z {\rm st}(\sigma Z, \alpha)] =  \sigma \Phi\big(\frac{\alpha-\mu}{\sigma}\big)  + \sigma \Phi\big(\frac{\alpha+\mu}{\sigma}\big)  \ .%=  \mathbb{E}_Z[Z {\rm relu}(\sigma Z + \mu - \alpha)] + \mathbb{E}_Z[Z {\rm relu}(\sigma Z - \mu - \alpha)] =  \sigma \Phi\big(\frac{\alpha-\mu}{\sigma}\big)  + \sigma \Phi\big(\frac{\alpha+\mu}{\sigma}\big)  \ .
\end{equation}
This implies that 
\begin{eqnarray}
    \mathbb{E}_{Z, \beta_0}\Big[\beta_0 \varphi\Big] &=&  \hat{w} \mathbb{E}_{Z}\Big[ \Phi\Big(\sqrt{\nu}\frac{\alpha\hat{\tau}-\hat{v}Z}{\hat{w}}\Big)  +  \Phi\Big(\sqrt{\nu}\frac{\alpha\hat{\tau}+\hat{v}Z}{\hat{w}}\Big)\Big] \nonumber \\
    &=& 2\hat{w} \mathbb{E}_{Z}\Big[ \Phi\Big(\sqrt{\nu}\frac{\alpha\hat{\tau}+\hat{v}Z}{\hat{w}}\Big)\Big] \ .
\end{eqnarray}
and also 
\begin{eqnarray}
    \mathbb{E}_{Z, \beta_0}\Big[Z \varphi\Big] &=& \hat{v}  \mathbb{E}_{ \beta_0}\Big[ \Phi\big(\frac{\alpha\hat{\tau}-\hat{w}\beta_0/\theta_0}{\hat{v}}\big) + \Phi\big(\frac{\alpha\hat{\tau}+\hat{w}\beta_0/\theta_0}{\hat{v}}\big) \Big] \nonumber \\
    &=& 2 \hat{v}  \mathbb{E}_{ \beta_0}\Big[\Phi\big(\frac{\alpha\hat{\tau}+\hat{w}\beta_0/\theta_0}{\hat{v}}\big) \Big]\ .
\end{eqnarray}
Next we use that 
\begin{equation}
    \mathbb{E}_X\Big[\Phi(aX + b)\Big] = \Phi\Big(\frac{b}{\sqrt{1 + a^2}}\Big)
\end{equation}
and get that 
\begin{equation}
    \mathbb{E}_{Z}\Big[ \Phi\Big(\sqrt{\nu}\frac{\alpha\hat{\tau}+\hat{v}Z}{\hat{w}}\Big)\Big] = \Phi\Big(\frac{\alpha\hat{\tau}}{\sqrt{\hat{w}^2/\nu + \hat{v}^2}}\Big)
\end{equation}
and 
\begin{equation}
    \mathbb{E}_{ \beta_0}\Big[\Phi\big(\frac{\alpha\hat{\tau}+\hat{w}\beta_0/\theta_0}{\hat{v}}\big) \Big] = \nu \Phi\Big(\frac{\alpha\hat{\tau}}{\sqrt{\hat{w}^2/\nu + \hat{v}^2}}\Big) + (1-\nu )\Phi\Big(\frac{\alpha\hat{\tau}}{\hat{v}}\Big) \ .
\end{equation}
Finally we use that
\begin{equation}
    \fl \mathbb{E}_Z[ {\rm relu}^2(\sigma Z + \mu - \alpha)]  = \Phi\big(\frac{\alpha-\mu}{\sigma}\big) \Big(\sigma^2 + (\mu-\alpha)^2 \Big) -  \sigma (\alpha-\mu) \frac{\rme^{-\frac{1}{2}\big(\frac{\alpha-\mu}{\sigma}\big) ^2}}{\sqrt{2\pi}}
\end{equation}
hence 
\begin{eqnarray}
     \fl \frac{1}{2}\mathbb{E}_{Z, \beta_0}\Big[\varphi^2\Big] &=&  \nu \bigg\{\Phi\big(\frac{\alpha\hat{\tau}}{\sqrt{\hat{v}^2 + \hat{w}^2/\nu}}\big) \Big(\hat{v}^2 + \hat{w}^2/\nu + \alpha^2\hat{\tau}^2 \Big) +\nonumber \\
    \fl &-&  \sqrt{\hat{v}^2 + \hat{w}^2/\nu} \alpha\hat{\tau} \frac{\rme^{-\frac{1}{2}\big(\frac{\alpha\hat{\tau}}{\sqrt{\hat{v}^2 + \hat{w}^2/\nu}}\big) ^2}}{\sqrt{2\pi}} \bigg\}+ \nonumber \\
    \fl &+& (1-\nu)\bigg\{\Phi\big(\frac{\alpha\hat{\tau}}{\hat{v}}\big) \Big(\hat{v}^2 + \alpha^2\hat{\tau}^2 \Big) -  \hat{v} \alpha \hat{\tau} \frac{\rme^{-\frac{1}{2}\big(\frac{\alpha\hat{\tau}}{\hat{v}}\big) ^2}}{\sqrt{2\pi}}\bigg\} \ .
\end{eqnarray}
If we introduce the short-hands 
\begin{equation}
    \chi_0 := \frac{\alpha\hat{\tau}}{\hat{v}}, \quad \chi_1 := \frac{\alpha\hat{\tau}}{\sqrt{\hat{v}^2 + \hat{w}^2/\nu}}
\end{equation}
then the first three RS equations take the more compact form 
\begin{eqnarray}
    w &=& 2\hat{w} \Phi(\chi_1)\\
    \tau &=& 2\hat{\tau}\Big\{ \nu \Phi(\chi_1) + (1-\nu) \Phi(\chi_0)\Big\}\\
    \frac{1}{2}(v^2 + w^2) &=&  \nu \{ (1 + 1/\chi_1^2) \Phi(\chi_1) - G(\chi_1)\} +\nonumber \\
    &+& (1-\nu) \{ (1 + 1/\chi_0^2) \Phi(\chi_0) - G(\chi_0)\}
\end{eqnarray}

\subsection{Explicit computation for the elastic net regularization}
\label{app_sec : computation_enet}
For the Elastic net regularization
\begin{equation}
    \prox_{\rmr(.)}(x, \hat{\tau}) = \frac{1}{1+\eta\hat{\tau}} {\rm st}(x, \alpha\hat{\tau})
\end{equation}
hence the replica symmetric equations of the main section can be easily recovered.

\section{Distributions}
\label{app_sec : distributions}
Here we show that the distributions 
\begin{eqnarray}
     \mathcal{P}_{\xi}(\Delta, t, h) &:=& \lim_{n\rightarrow \infty }\mathbb{E}_{\mathcal{D}}\Big[\frac{1}{n}\sum_{i=1}^n \delta_{\Delta, \Delta_i}\delta(t - T_i)\delta (h - \mathbf{X}_i'\hat{\bbeta})\Big]\ ,   \\
     \mathcal{P}_{\varphi}(x) &:=& \lim_{n\rightarrow \infty }\mathbb{E}_{\mathcal{D}}\Big[\frac{1}{p}\sum_{k=1}^p \delta (x - \mathbf{e}_{k}'\hat{\bbeta})\Big]
\end{eqnarray}
indeed satisfy 
\begin{eqnarray}
\label{app : law_xi_rs}
    \mathcal{P}_{\xi}(\Delta, t, h) &=& \mathbb{E}_{\Delta', T', Z'_0, Q'} \Big[\delta(t - T')\delta_{\Delta, \Delta'}\delta\big(h - \xi_{\star}\big)\Big]\\
\label{app : law_varphi_rs}
    \mathcal{P}_{\varphi}(x) &=& \mathbb{E}_{\beta_0, Z} \Big[\delta\big(x - \varphi_{\star}\big)\Big], 
\end{eqnarray}
as stated in the main text (\ref{def : law_xi_rs}, \ref{def : law_varphi_rs}). 
For the identity (\ref{app : law_xi_rs}), observe that $ \mathcal{P}_{\xi}(\Delta, t, h)$ is the functional order parameter of the replica theory (\ref{def : functional_order_parameter}), hence, at the replica symmetric saddle point we have (\ref{def : rs_P_gamma_r}) and after taking the limit $r\rightarrow 0$ and $\gamma\rightarrow\infty$ we obtain (\ref{def : P_rs}), which coincides with (\ref{app : law_xi_rs}) as 
\begin{eqnarray}
     \mathcal{P}_{\xi}(\Delta,t,h)  &=& \mathbb{E}_{Z_0,Z}\Big[f(\Delta,t|Z_0)\delta\Big(h - \xi_{\star}\Big)\Big]  =\nonumber \\
     &=& \mathbb{E}_{\Delta', T', Z'_0,Z'}\Big[\delta_{\Delta,\Delta'}\delta(t-T')\delta\Big(h - \xi_{\star}\Big)\Big] \ .
\end{eqnarray}
The identity (\ref{app : law_varphi_rs}) requires more care. Via Laplace integration we may write 
\begin{equation}
    \frac{1}{p}\sum_{k=1}^p \delta (x - \mathbf{e}_{k}'\hat{\bbeta}) =\lim_{\gamma\rightarrow\infty} \frac{1}{p}\sum_{i=1}^p \sum_{k=1}^p\frac{\int \rme^{-\gamma \mathcal{H}(\bbeta)} \ \delta (x - \mathbf{e}_{\mu}'\bbeta) \ \rmd \bbeta}{\int \rme^{-\gamma \mathcal{H}(\bbeta)} \ \rmd \bbeta} \ .
\end{equation}
To compute the expectation over the data-set we use the following alternative replica identity
\begin{equation}
     \fl \mathcal{P}_{\varphi}(x)  =\lim_{r\rightarrow 0}\lim_{\gamma\rightarrow\infty}\lim_{p\rightarrow\infty}\frac{1}{ r \ p}\sum_{k=1}^p\mathbb{E}_{\mathcal{D}}\bigg[\int \rme^{-\gamma\sum_{\alpha=1}^r \mathcal{H}(\bbeta_{\alpha}|\mathcal{D})} \ \delta (x - \mathbf{e}_{k}'\bbeta_{1}) \ \prod_{\alpha=1}^r\rmd \bbeta_{\alpha} \bigg]\ .
\end{equation}
Inserting the functional delta measure as in (\ref{def : functional_delta}) and taking the expectation over $\mathcal{D}$, i.e. the data-set,  we obtain 
\begin{equation}
    \hspace{-1cm} \mathcal{P}_{\varphi}(x) =  \int \rme^{ -n \psi\big[\{\mathcal{P}_{\alpha},\hat{\mathcal{P}}_{\alpha}\}_{\alpha=1}^r,\hat{\mathbf{C}},\mathbf{C}\big] } \mathcal{J}(\gamma, \hat{\mathbf{C}})\prod_{\alpha=1}^r  \ \mathcal{D}\hat{\mathcal{P}}_{\alpha} \ \mathcal{D}\mathcal{P}_{\alpha} \rmd \mathbf{C} \rmd \hat{\mathbf{C}}
\end{equation}
where $\psi$ is defined as in (\ref{def : psi}) and 
\begin{equation}
    \hspace{-1cm} \mathcal{J}(\gamma, \hat{\mathbf{C}}) :=  \frac{1}{ p}\sum_{k=1}^p\frac{\int \rme^{-\rmi p{\rm Tr}( \hat{\mathbf{C}}\mathbf{C}(\{\bbeta_{\alpha}\})) -\gamma \rmr(\bbeta_{\alpha})}\delta(x - \mathbf{e}_{k}'\bbeta_1)\prod_{\alpha=1}^r \rmd \bbeta_{\alpha}}{\int \rme^{-\rmi p{\rm Tr}( \hat{\mathbf{C}}\mathbf{C}(\{\bbeta_{\alpha}\})) -\gamma \rmr(\bbeta_{\alpha})}\prod_{\alpha=1}^r \rmd \bbeta_{\alpha}} \ .
\end{equation}
Now taking $i\hat{\mathbf{C}} = \frac{1}{2}\mathbf{D}$, and replica symmetry, we notice that the expression above can be simplified, at the saddle point, because the regularizer is separable, i.e. $r(\bx) = \sum_{l=1}^p r(x_l)$, 
\begin{eqnarray}
    \fl &&\mathcal{J}_{RS}^{(r)}(\gamma, \hat{m}, \hat{q}, \hat{\rho}) =  \nonumber \\
   \fl  &&\frac{1}{ p}\sum_{k=1}^p\frac{ \int \rme^{-\sum_{\rho=1}^r \big( \hat{m}\beta_{0,k}\beta_{\rho, k}  + \frac{1}{2} (\hat{\rho} +\hat{q})  \beta_{\rho, k}^2\big) + \frac{1}{2}\hat{q}\big(\sum_{\rho=1}^r  \bbeta_{\rho, k}\big)^2 -\gamma \rmr(\beta_{\alpha,k })}\delta(x - \beta_{1,k})\prod_{\alpha=1}^r \rmd \beta_{\alpha, k} }{\int \rme^{-\sum_{\rho=1}^r \big( \hat{m}\beta_{0,k}\beta_{\rho, k}  + \frac{1}{2} (\hat{\rho} +\hat{q})  \beta_{\rho, k}^2\big) + \frac{1}{2}\hat{q}\big(\sum_{\rho=1}^r  \bbeta_{\rho, k}\big)^2 -\gamma \rmr(\beta_{\alpha,k})}\prod_{\alpha=1}^r \rmd \beta_{\alpha, k}} \ . \nonumber
\end{eqnarray}
At this point Gaussian linearization gives
\begin{eqnarray}
    \fl &&\mathcal{J}_{RS}^{(r)}(\gamma, \hat{m}, \hat{q}, \hat{\rho}) =  \mathbb{E}_{\beta_0}\bigg[\mathbb{E}_{Z}\bigg[ \int \rme^{- \big( \{\hat{m}\beta_{0} + \sqrt{\hat{q}} Z \}\beta  + \frac{1}{2} (\hat{\rho} +\hat{q}) \beta^2\big)  -\gamma \rmr(\beta)}\delta(x - \beta)\rmd \beta \times   \\
    \fl &&\hspace{0.5cm}\times\bigg(\int \rme^{- \big( \{\hat{m}\beta_{0,k} + \sqrt{\hat{q}} Z \}\beta  + \frac{1}{2} (\hat{\rho} +\hat{q}) \beta^2\big)  -\gamma \rmr(\beta)}\rmd \beta\bigg)^{r-1}\bigg] \times \\
    \fl &&\hspace{0.5cm}\times\mathbb{E}_{Z}\bigg[ \bigg(\int \rme^{- \big( \{\hat{m}\beta_{0} + \sqrt{\hat{q}} Z \}\beta  + \frac{1}{2} (\hat{\rho} +\hat{q}) \beta^2\big)  -\gamma \rmr(\beta)}\rmd \beta\bigg)^{r}\bigg]^{-1} \bigg]\ .
\end{eqnarray}
Hence, after taking the limit $r\rightarrow 0$, we get 
\begin{equation}
    \fl \mathcal{J}_{RS}(\gamma, \hat{m}, \hat{q}, \hat{\rho}) :=  \mathbb{E}_{Z, \beta_0}\bigg[ \frac{\int \rme^{- \big( \{\hat{m}\beta_{0} + \sqrt{\hat{q}} Z \}\beta  + \frac{1}{2} (\hat{\rho} +\hat{q}) \beta^2\big)  -\gamma \rmr(\beta)}\delta(x - \beta)\rmd \beta}{\int \rme^{- \big( \{\hat{m}\beta_{0,} + \sqrt{\hat{q}} Z \}\beta  + \frac{1}{2} (\hat{\rho} +\hat{q}) \beta^2\big)  -\gamma \rmr(\beta)}\rmd \beta}\bigg]\ .
\end{equation}
After taking the rescaling
\begin{equation}
    1/\hat{\tau} = (\hat{\rho} +\hat{q}) /\gamma, \quad \hat{m} = \gamma \hat{m}, \quad \hat{q} = \gamma^2 \hat{q}
\end{equation}
and the limit $\gamma\rightarrow\infty$, we get 
\begin{equation}
   \mathcal{P}_{\varphi}(x) =\mathbb{E}{Z, \beta_0}\Big[\delta \Big(x - {\rm prox}_{\rmr(.)} \big(\hat{m}_{\star}\beta_0 +  \sqrt{\hat{q}_{\star}} Z, \hat{\tau}_{\star}\big)\Big) \Big]\ .
\end{equation}
The expression above reduces to (\ref{app : law_varphi_rs}) after the change of variables
\begin{equation}
     \hat{w} = -\hat{\tau} \hat{m} \theta_0, \quad \hat{v} =  \hat{\tau}\sqrt{\hat{q}}, 
\end{equation}
and using the definition (\ref{def : varphi}).

\section{Coordinate Wise Minimization algorithm for pathwise solution}
\label{app : cd_algorithm}
Consider the function
\begin{equation}
\label{app : objective}
    g(\bbeta,\lambda) := \frac{1}{n} \sum_{i=1}^n \Big\{\Lambda(T_i)\rme^{\mathbf{X}_i'\bbeta} - \Delta_i\mathbf{X}_i'\bbeta  -\Delta_i \log \lambda(T_i)\Big\}
\end{equation}
Suppose we want to minimize the objective function
\begin{equation}
    f(\bbeta, \lambda) := g(\bbeta, \lambda) +\rmr(\bbeta) 
\end{equation}
with $\rmr$ a separable convex regularization function. In our case we will be interested in 
\begin{equation}
\label{app : enet}
     \rmr(\bbeta)  := \alpha |\bbeta| + \frac{1}{2}\eta \|\bbeta\|^2 \ .
\end{equation}
This can be done via coordinate descent. First we minimize over $\lambda$, obtaining the Nelson-Aalen estimator 
\begin{equation}
\label{app: NA_est}
    \hat{\Lambda}_n(t) =  {\rm NA}(\{T_j, \mathbf{X}_j'\bbeta\} ) := \sum_{i=1}^n \frac{\Delta_i \Theta(t-T_i)}{\sum_{j=1}^n \Theta(T_j-T_i)\rme^{\mathbf{X}_j'\bbeta}}   \ .
\end{equation}
Then we minimize over $\bbeta$ and so on an so forth, until a fixed point, i.e. the updated values for $\bbeta^{t+1}$ and $\Lambda^{t+1}(T_1), \dots, \Lambda^{t+1}(T_n)$ are within a user defined tolerance from their previous value.
To compute the minimizer in $\bbeta$ at fixed $\Lambda$, we use the algorithm proposed in \cite{Zou_05}.
The idea is to reduce the problem to an iterative regularized least squared regression.
Let us define 
\begin{equation}
\label{app : objective_beta}
    \ell(\bbeta, \Lambda) := \frac{1}{n} \sum_{i=1}^n \Big\{\Lambda(T_i)\rme^{\mathbf{X}_i'\bbeta} - \Delta_i\mathbf{X}_i'\bbeta \Big\} \ .
\end{equation}
A second order expansion in $\bxi$ around $\bphi$ gives the following approximation for $\ell$
\begin{equation}
    \tilde{\ell}(\bbeta, \Lambda, \bphi) = \ell(\bphi, \Lambda) + \mathbf{s}(\bphi,\Lambda)' (\bbeta - \bphi) + \frac{1}{2}  (\bbeta- \bphi)'\bm{H}(\bphi, \Lambda) (\bbeta - \bphi)
\end{equation}
where 
\begin{eqnarray}
    \mathbf{s}(\bphi,\Lambda) &=&  \frac{1}{n} \sum_{i=1}^n \Big\{\Lambda(T_i)\rme^{\mathbf{X}_i'\bphi} - \Delta_i\Big\}\mathbf{X}_i = \frac{1}{n} (\mathbf{W} - {\rm diag}(\bm{\Delta})\mathbf{X} \\
    \mathbf{M}(\bphi,\Lambda) &=&  \frac{1}{n} \sum_{i=1}^n \Lambda(T_i)\rme^{\mathbf{X}_i'\bphi} \mathbf{X}_i\mathbf{X}_i' = \frac{1}{n} \mathbf{X}'\mathbf{W}\mathbf{X}
\end{eqnarray}
with $\mathbf{W} := {\rm diag}\big(\Lambda(T_i)\exp\{\mathbf{X}_i'\bphi\}\big)$.
At this points we use a coordinate descent strategy to solve the regularized least square problem 
\begin{equation}
    \bbeta^{t+1} = \underset{\bvarphi}{\arg\min} \Big\{f(\bvarphi) \Big\}, \quad f(\bvarphi) = \tilde{\ell}(\bvarphi, \Lambda^{t}, \bbeta^{t}) + r(\bvarphi), 
\end{equation} 
i.e. we minimize along each component keeping the remaining components fixed
\begin{eqnarray}
   \fl && \frac{\partial}{\partial \varphi_k}f(\bvarphi)  = \nonumber\\
   \fl && s_k(\bbeta^{t}, \Lambda^{t}) + \varphi_k \ \mathbf{e}_k'\mathbf{M}(\bbeta^{t}, \Lambda^{t})\mathbf{e}_k+ \mathbf{e}_k'\mathbf{M}(\bbeta^{t}, \Lambda^{t})\big\{(\bm{I} - \mathbf{e}_k\mathbf{e}_k')\bvarphi -\bbeta^{t}\big\} + r'(\varphi_k) = 0 
\end{eqnarray}
which is solved by 
\begin{eqnarray}
    \fl &&\varphi_k^{t+1} = \nonumber\\
    \fl && {\rm prox}_{r(.)}\Big( \frac{\mathbf{e}_k'\mathbf{M}(\bbeta^{t}, \Lambda^{t})\big\{\bbeta^{t} -(\bm{I} - \mathbf{e}_k\mathbf{e}_k')\bvarphi^{t}\big\}- s_k(\bbeta^{t}, \Lambda^{t})}{\mathbf{e}_k'\mathbf{M}(\bbeta^{t}, \Lambda^{t})\mathbf{e}_k}, \frac{1}{\mathbf{e}_k'\mathbf{M}(\bbeta^{t}, \Lambda^{t})\mathbf{e}_k}  \Big) \ .
\end{eqnarray}
Notice that in the case of the elastic net penalization we get 
\begin{equation}
    \varphi_k^{t+1} =\frac{1}{1 + \eta \tau_k} {\rm st}\big( \psi_k, \alpha  \tau_k\big)
\end{equation}
where 
\begin{eqnarray}
    \fl \psi_k &:=& \big\{\mathbf{e}_k'\mathbf{M}(\bbeta^{t}, \Lambda^{t})\big\{\bbeta^{t} -(\bm{I} - \mathbf{e}_k\mathbf{e}_k')\bvarphi^{t}\big\}- s_k(\bbeta^{t}, \Lambda^{t})\big\}/\mathbf{e}_k'\mathbf{M}(\bbeta^{t}, \Lambda^{t})\mathbf{e}_k \\
    \fl 1/\tau_k &:=&\mathbf{e}_k'\mathbf{M}(\bbeta^{t}, \Lambda^{t})\mathbf{e}_k
\end{eqnarray}
and ${\rm st}$ is the soft thresholding operator.

\section{Derivation of GAMP via Belief Propagation}

\label{app : amp_cox}
Consider the following optimization problem 
\begin{equation}
    \hat{\bbeta} = \underset{\bbeta}{\arg\min} \Big\{\sum_{i=1}^n g\big(\mathbf{X}_i'\bbeta,Y_i\big) + \sum_{\mu=1}^p r(\beta_{\mu})\Big\}
\end{equation}
where $\mathbf{X}_i\sim\mathcal{N}\big(0,\frac{1}{p}\bm{1}\big)$ and $Y_i|\mathbf{X}_i \sim f(.|\mathbf{X}_i'\bbeta_0)$.
Define the Gibbs measure
\begin{equation}
    f_{\gamma}(\bbeta|\{Y_i,\mathbf{X}_i\})  =  \frac{\rme^{-\gamma  \big\{\sum_{i=1}^n g(\mathbf{X}_i'\bbeta,Y_i) + \lambda \sum_{\mu=1}^p r(\bbeta_{\mu})\big\}}}{Z_{\gamma}(\{Y_i,\mathbf{X}_i\})} \ . 
\end{equation}
Such a distribution corresponds to a fully connected bypartite factor graph, with $p$ variable nodes and $n$ factor nodes.
The (loopy) belief propagation equations read
\begin{eqnarray}
    m_{\mu \rightarrow i}(\bbeta_{\mu}) &\propto& \rme^{-\gamma  r(\bbeta_{\mu})} \prod_{j\neq i} m_{j\rightarrow \mu}(\beta_{\mu})\\
    m_{i\rightarrow \mu}(\bbeta_{\mu}) &\propto& \int  \rme^{-\gamma g(\mathbf{X}_i'\bbeta,Y_i)} \prod_{\alpha \neq \mu} m_{\alpha \rightarrow i}(\beta_{\alpha}) \rmd \beta_{\alpha} \ .
\end{eqnarray}

\subsection{The tempered Moreau envelope and proximal operator}
Let us introduce the tempered Moreau envelope as 
\begin{equation}
    \mathcal{M}_{f(.)} (x, \alpha, \gamma) := -\frac{1}{\gamma} \log \int \rme^{-\gamma\big\{\frac{1}{2}\frac{(z-x)^2}{\alpha} + f(z) \big\}} \rmd z \ .
\end{equation}
Notice that 
\begin{eqnarray}
    \dot{\mathcal{M}}_{f(.)} (x, \alpha, \gamma) = \frac{1}{\alpha}\big( x -\mathbb{E}_{\gamma}[Z]\big)\\
    \ddot{\mathcal{M}}_{f(.)} (x, \alpha, \gamma) = \frac{1}{\alpha} - \frac{\gamma}{\alpha^2}\mathbb{V}_{\gamma}[Z]\ .
\end{eqnarray}
In analogy with the proximal operator, we define 
\begin{equation}
    \fl \prox_{f(.)}(x, \alpha, \gamma) := \mathbb{E}_{\gamma}[Z] = \frac{\int z \rme^{-\gamma\big\{\frac{1}{2}\frac{(z-x)^2}{\alpha} + f(z) \big\}}\rmd z}{\int \rme^{-\gamma\big\{\frac{1}{2}\frac{(z-x)^2}{\alpha} + f(z) \big\}}\rmd z } = x - \alpha \dot{\mathcal{M}}_{f(.)} (x, \alpha, \gamma) \ .
\end{equation}
Based on the previously derived relationships we notice that 
\begin{equation}
    \prox'_{f(.)}(x, \alpha, \gamma) = \mathbb{V}_{\gamma}[Z] \gamma / \alpha  \ .
\end{equation}
In the limit $\gamma \rightarrow\infty$ we have 
\begin{equation}
    \lim_{\gamma\rightarrow\infty}\mathcal{M}_{f(.)} (x, \alpha, \gamma) := \mathcal{M}_{f(.)} (x, \alpha) 
\end{equation}
furthermore 
\begin{eqnarray}
     \fl \lim_{\gamma\rightarrow\infty} \dot{\mathcal{M}}_{f(.)} (x, \alpha, \gamma) = \frac{1}{\alpha}\big( x -\prox_{f(.)}(x, \alpha)\big) = \dot{f}\big(\prox_{f(.)}(x, \alpha)\big)\\
    \fl \lim_{\gamma\rightarrow\infty}\ddot{\mathcal{M}}_{f(.)} (x, \alpha, \gamma) = \frac{\ddot{f}\big(\prox_{f(.)}(x, \alpha)\big)}{1 + \alpha \ddot{f}\big(\prox_{f(.)}(x, \alpha)\big)}\ .
\end{eqnarray}
The last identity is obtained by differentiating $\dot{f}\big(\prox_{f(.)}(x, \alpha)\big)$ with respect to $x$.
Notice that this implies 
\begin{equation}
     \lim_{\gamma\rightarrow\infty}\gamma \mathbb{V}_{\gamma}[Z] = \alpha \frac{1}{1 + \alpha \ddot{f}\big(\prox_{f(.)}(x, \alpha)\big)} \ .
\end{equation}

\subsection{Approximation of the messages}
We start by re-writing the self consistent condition of the messages as 
\begin{eqnarray}
    \fl m_{i\rightarrow \mu}(\bbeta_{\mu}) &\propto& \int  \rme^{-\gamma g \big(X_{i,\mu}\beta_{\mu} + \sum_{\alpha\neq \mu} X_{i,\alpha}\beta_{\alpha},Y_i\big)} \prod_{\alpha \neq \mu} \rmd m_{\alpha\rightarrow i}(\beta_{\alpha}) \nonumber \\
    \fl  &\propto& \int  \rme^{i\hat{\xi}\big(\xi- X_{i,\mu}\beta_{\mu}\big) - \sum_{\alpha\neq \mu} \psi_{\beta_{\alpha}}(i\hat{\xi}X_{i,\alpha})  - \gamma g(\xi,Y_i) }\ \rmd \xi \rmd \hat{\xi}  
\end{eqnarray}
where 
\begin{equation}
    \psi_{X}(q) := \log \mathbb{E}_X\Big[\rme^{qX}\Big]
\end{equation}
is the cumulant generating function.
Taking a quadratic approximation of $ \psi_{\beta_{\alpha}}$ gives
\begin{eqnarray}
    \fl m_{i\rightarrow \mu}(\bbeta_{\mu})&\approx & \int  \rme^{i\hat{\xi}\big(\xi-  X_{i,\mu}\beta_{\mu} - \sum_{\alpha\neq \mu} X_{i,\alpha}\hat{\beta}_{\alpha\rightarrow i}\big)  - \gamma  g(\xi,Y_i) - \frac{1}{2} \frac{1}{\gamma}\big(\sum_{\alpha\neq \mu} \nu_{\alpha \rightarrow i} X_{i,\alpha}^2\big) \hat{\xi}^2}   \ \rmd \xi \rmd \hat{\xi}= \nonumber  \\
    \fl  &\propto& \int  \rme^{-\gamma\Big\{\frac{1}{2}\frac{\big(\xi-  X_{i,\mu}\beta_{\mu} - \sum_{\alpha\neq \mu} X_{i,\alpha}\hat{\beta}_{\alpha\rightarrow i}\big)^2}{\sum_{\alpha\neq \mu}\nu_{\alpha\rightarrow i} X_{i,\alpha}^2}  +  g(\xi,Y_i)  \Big\}}   \ \rmd \xi  \ ,
\end{eqnarray}
where $\hat{\beta}_{\alpha\rightarrow i}, \tau_{\alpha\rightarrow i} $ are respectively the mean and the variance (rescaled by the inverse temperature $\gamma$) of the belief $m_{\alpha \rightarrow i}$. 
We notice that 
\begin{equation}
    m_{i\rightarrow \mu}(\bbeta_{\mu}) \approx \rme^{-\gamma \dot{\mathcal{M}}_{f(.)} (X_{i,\mu}(\beta_{\mu} - \hat{\beta}_{\mu\rightarrow i}) +\hat{\xi}_i, \tau_i - \nu_{\mu\rightarrow i} X_{i,\mu}^2, \gamma)}
\end{equation}
where 
\begin{eqnarray}
    \xi_{i} &=& \sum_{\alpha} X_{i,\alpha}\hat{\beta}_{\alpha\rightarrow i}\\
    \tau_{i} &=&  \sum_{\alpha} \nu_{\alpha\rightarrow i} X_{i,\alpha}^2  \ .
\end{eqnarray}
Assuming $  X_{i,\mu}(\beta_{\mu} - \hat{\beta}_{\mu\rightarrow i})$ to be small, we take a quadratic approximation of $\mathcal{M}_{g(., T_i)}(., \tau_i - \nu_{\mu\rightarrow i} X_{i,\mu}^2, \gamma )$ around $\xi_{i}$ obtaining
\begin{eqnarray}
    \fl &&\mathcal{M}_{g(., T_i)} (X_{i,\mu}(\beta_{\mu} - \hat{\beta}_{\mu\rightarrow i}) +\xi_i, \tau_i - \nu_{\mu\rightarrow i} X_{i,\mu}^2, \gamma) \approx \nonumber \\
    \fl &&\hspace{1cm}\mathcal{M}_{g(., T_i)}( \xi, \tau_i - \nu_{\mu\rightarrow i} X_{i,\mu}^2, \gamma) + \nonumber\\
    \fl &&\hspace{1cm} \dot{\phi}_{i\rightarrow \mu}X_{i,\mu}(\beta_{\mu}  - \hat{\beta}_{\mu\rightarrow i}) +\frac{1}{2} \ddot{\phi}_{i\rightarrow\mu} X^2_{i,\mu}(\beta_{\mu} - \hat{\beta}_{\mu\rightarrow i})^2
\end{eqnarray}
where 
\begin{eqnarray}
     \dot{\phi}_{i\rightarrow \mu} &=& \dot{\mathcal{M}}_{g(., T_i)} (\xi_i, \tau_i - \nu_{\mu\rightarrow i} X_{i,\mu}^2, \gamma)\\
     \ddot{\phi}_{i\rightarrow \mu} &=& \ddot{\mathcal{M}}_{g(., T_i)} (\xi_i, \tau_i - \nu_{\mu\rightarrow i} X_{i,\mu}^2, \gamma) \ .
\end{eqnarray}
Then 
\begin{equation}
    m_{i\rightarrow \mu}(\bbeta_{\mu}) \approx \exp\big\{-\gamma \Big( \frac{1}{2}\hat{\tau}_{i \rightarrow \mu }\beta_{\mu}^2 - s_{i \rightarrow \mu }\beta_{\mu} \Big)\big\} \ .%= \mathcal{N}\big(\hat{\beta}_{i\rightarrow\mu}, \nu_{i\rightarrow \mu})
\end{equation}
where 
\begin{eqnarray}
    s_{i \rightarrow \mu } &=& -(\dot{\phi}_{i\rightarrow \mu}X_{i,\mu} -\ddot{\phi}_{i\rightarrow\mu} X^2_{i,\mu}\hat{\beta}_{\mu\rightarrow i})\\
    1 / \hat{\tau}_{i \rightarrow \mu} &=& \ddot{\phi}_{i\rightarrow\mu} X^2_{i,\mu} \ .
\end{eqnarray}
%with 
%\begin{eqnarray}
%    \hat{\beta}_{i\rightarrow\mu} &=& \hat{\beta}_{\mu\rightarrow i} -  \dot{\phi}_{i\rightarrow 5\mu}X_{i,\mu} / \ddot{\phi}_{i\rightarrow\mu} X^2_{i,\mu}  \\
%    \nu_{i\rightarrow \mu} &=&  1.0  / \big(\gamma \ddot{\phi}_{i\rightarrow\mu} X^2_{i,\mu}\big)  \ ,
%\end{eqnarray}
As a consequence
\begin{equation}
    m_{\mu\rightarrow i}(\beta_{\mu}) \approx \exp\Big\{-\gamma \Big(  r(\beta_{\mu})  + \frac{1}{2}   \beta_{\mu}^2 / \hat{\tau}_{\mu\rightarrow i} -  s_{\mu\rightarrow i}\beta_{\mu}\Big)\Big\}  \ .
\end{equation}
where 
\begin{eqnarray}
    s_{\mu\rightarrow i} &=& \sum_{j\neq i}s_{j \rightarrow \mu }\\
    1 / \hat{\tau}_{\mu\rightarrow i} &=&\sum_{j\neq i }\hat{\tau}^{-1}_{j \rightarrow \mu } \ .
\end{eqnarray}
Finally 
\begin{eqnarray}
    \hat{\beta}_{\mu\rightarrow i} = \hat{\tau}_{\mu\rightarrow i } s_{\mu\rightarrow i}  -  \hat{\tau}_{\mu\rightarrow i } \dot{\mathcal{M}}_{r(.)} \big( s_{\mu\rightarrow i} \hat{\tau}_{\mu\rightarrow i },  \hat{\tau}_{\mu\rightarrow i } , \gamma\big)  \\
    \nu_{\mu\rightarrow i} = \hat{\tau}_{\mu\rightarrow i } \Big\{ 1 -   \hat{\tau}_{\mu\rightarrow i }  \ddot{\mathcal{M}}_{r(.)} \big( s_{\mu\rightarrow i} \hat{\tau}_{\mu\rightarrow i },  \hat{\tau}_{\mu\rightarrow i } , \gamma\big)\Big\} \ ,
\end{eqnarray}
and as a consequence 
\begin{eqnarray}
    \fl \hat{\beta}_{\mu\rightarrow i} = \prox_{r(.)} \big( s_{\mu\rightarrow i} \hat{\tau}_{\mu\rightarrow i },  \hat{\tau}_{\mu\rightarrow i } , \gamma\big)  = \prox_{r(.)} \big( \psi_{\mu\rightarrow i},  \hat{\tau}_{\mu\rightarrow i } , \gamma\big)  \\
   \fl  \nu_{\mu\rightarrow i} = \hat{\tau}_{\mu\rightarrow i } \prox'_{r(.)} \big( s_{\mu\rightarrow i} \hat{\tau}_{\mu\rightarrow i },  \hat{\tau}_{\mu\rightarrow i } , \gamma\big) = \hat{\tau}_{\mu\rightarrow i } \prox'_{r(.)} \big( \psi_{\mu\rightarrow i},  \hat{\tau}_{\mu\rightarrow i } , \gamma\big) \ , 
\end{eqnarray}
where 
\begin{equation}
     \psi_{\mu\rightarrow i} := s_{\mu\rightarrow i} \hat{\tau}_{\mu\rightarrow i } = \hat{\beta}_{\mu \rightarrow i}  - \hat{\tau}_{\mu\rightarrow i }X_{i,\mu}\dot{\phi}_{i\rightarrow\mu} \ .
\end{equation}

\subsection{TAP approximation}

First of all we neglect all the corrections $\nu_{\mu\rightarrow i} X_{i,\mu}^2$ so that 
\begin{eqnarray}
     \dot{\phi}_{i\rightarrow \mu} &\approx& \dot{\phi}_{i} = \dot{\mathcal{M}}_{g(., T_i)} (\xi_i, \tau_i, \gamma)\\
     \ddot{\phi}_{i\rightarrow \mu} &\approx& \ddot{\phi}_{i} =  \ddot{\mathcal{M}}_{g(., T_i)} (\xi_i, \tau_i, \gamma) \ .
\end{eqnarray}
Then
\begin{eqnarray}
    s_{i \rightarrow \mu } &\approx& -(\dot{\phi}_{i}X_{i,\mu} -\ddot{\phi}_{i} X^2_{i,\mu}\hat{\beta}_{\mu\rightarrow i})\\
    \hat{\tau}_{i \rightarrow \mu} &\approx& \ddot{\phi}_{i} X^2_{i,\mu}  \ .
\end{eqnarray}
Notice that 
\begin{eqnarray}
    s_{\mu\rightarrow i} &=& \sum_{j}s_{j \rightarrow \mu } - s_{i\rightarrow \mu} = s_{\mu}- s_{i\rightarrow \mu}\\
    1/\hat{\tau}_{\mu\rightarrow i} &=&\sum_{j}1/(\hat{\tau}_{j \rightarrow \mu } - \hat{\tau}_{i \rightarrow \mu }) \approx  1 / \hat{\tau}_{\mu}  \ , 
\end{eqnarray}
hence 
\begin{eqnarray}
   \fl  \hat{\beta}_{\mu\rightarrow i} &\approx& \prox_{\rmr(.)}\big( \psi_{\mu} , \hat{\tau}_{\mu} , \gamma\big)  -  \prox'_{\rmr(.)}\big( \psi_{\mu} , \hat{\tau}_{\mu} , \gamma\big)s_{i\rightarrow\mu} \hat{\tau}_{\mu}  = \hat{\beta}_{\mu}  -  \nu_{\mu} s_{i\rightarrow\mu}\\
    \fl \nu_{\mu\rightarrow i} &\approx& \nu_{\mu} = \prox'_{\rmr(.)}\big( \psi_{\mu} , \hat{\tau}_{\mu} , \gamma\big)\hat{\tau}_{\mu}\  , 
\end{eqnarray}
where 
\begin{equation}
     \psi_{\mu} := s_{\mu} \hat{\tau}_{\mu} = \hat{\beta}_{\mu }  - \hat{\tau}_{\mu }X_{i,\mu}\dot{\phi}_{i} 
\end{equation}
and 
\begin{equation}
    \hat{\beta}_{\mu} = \prox_{\rmr(.)}\big( \psi_{\mu} , \hat{\tau}_{\mu} , \gamma\big) \ .
\end{equation}
Finally
\begin{eqnarray}
    \tau_{i} &\approx&  \sum_{\alpha} \nu_{\alpha} X_{i,\alpha}^2 \nonumber \\
    \xi_{i} &\approx& \sum_{\alpha} X_{i,\alpha}\hat{\beta}_{\alpha} + \sum_{\alpha} X_{i,\alpha}^2 \nu_{\alpha} \dot{\phi}_i %= \prox_{g(., T)} ( \sum_{\alpha} X_{i,\alpha}\hat{\beta}_{\alpha}, \tau_i, \gamma ) \ .
\end{eqnarray}
\newpage
\subsection{TAP algorithm for GLMs}
The TAP  algorithm can be schematized as follows
\begin{algorithm}
\caption{TAP algorithm for Bayes Estimator in GLMs}\label{ alg : tap_GLM}
\begin{algorithmic}
\Require{$\hat{\bbeta}^{0} \gets \bm{0}_p$,  tol $\gets 1.0e-8$, max epochs $\gets 300$}
\State $\bm{\tau}^0 = \bm{0}_n, \ \hat{\bm{\tau}}^0 = \bm{0}_p, \bm{\xi}^{0} = \mathbf{X}\hat{\bbeta}^{0}$
\State flag $=$ True,  err $= \infty$, t $= 0$ 
\While{err $\geq$ tol and flag}
\State $t \gets t + 1$
\State $\xi_{i}^{t} \gets \mathbf{X}_{i}'\hat{\bbeta}^{t-1} + \tau_i \dot{\mathcal{M}}_{g(. , T)}(\xi^{t-1}, \tau^{t-1}_i, \gamma),\ \forall i$
\State $1/\hat{\tau}_{\mu}^{t} \gets \sum_{i=1}^n X^2_{i,\mu} \ddot{\mathcal{M}}_{g(. , T)}(\xi^t_i, \tau^{t-1}_i, \gamma), \ \forall \mu $ %\frac{\ddot{g}(\xi_i^{t}, T)}{1 + \tau_i^{t} \ddot{g}(\xi_i^{t}, T)}
\State $\psi_{\mu}^t \gets  \hat{\beta}_{\mu}^{t-1} - \hat{\tau}_{\mu}^t\sum_{i=1}^n \dot{\mathcal{M}}_{g(. , T)}(\xi^t_i, \tau^{t-1}_i, \gamma)X_{i,\mu}, \ \forall \mu$
\State $\hat{\beta}_{\mu}^{t}\gets \prox_{r(.)}(\psi_{\mu}^t, \hat{\tau}_{\mu}^t, \gamma), \ \forall \mu $
\State $\tau_i^{t} \gets \sum_{\alpha = 1}^p   \hat{\tau}_{\alpha}^t \ \prox'_{r(.)}(\psi_{\alpha}^{t}, \hat{\tau}_{\alpha}^t) X_{i,\alpha}^2, \ \forall i$
\State err $\gets \sqrt{ \|\hat{\bbeta}^{t} - \hat{\bbeta}^{t-1}\|_2^2 + \|\bm{\tau}^t - \bm{\tau}^{t-1}\|_2^2 +  \|\bm{\xi}^t - \bm{\xi}^{t-1}\|_2^2 + \|\hat{\bm{\tau}}^t - \hat{\bm{\tau}}^{t-1}\|_2^2}$
\If{t $\geq$ max epochs}
\State{flag $=$ False}
\EndIf
\EndWhile
\end{algorithmic}
\end{algorithm}

The algorithm above is, in general, difficult to implement, since it requires to compute $\ddot{\mathcal{M}}_{g(. , T)}(., a, \gamma)$, $\dot{\mathcal{M}}_{g(. , T)}(., a, \gamma), \prox_{r(.)}(., a, \gamma)$, $ \prox'_{r(.)}(., a, \gamma)$ at finite $\gamma<\infty$. These quantities are not directly available (an exception is the linear model with ridge regularization) and must be approximated via a Montecarlo scheme, or by numerical integration.
In the zero temperature limit, however, the algorithm simplify to 
\begin{algorithm}
\caption{TAP algorithm for MAP in GLMs}\label{alg : tap_MAP_GLM}
\begin{algorithmic}
\Require{$\hat{\bbeta}^{0} \gets \bm{0}_p$,  tol $\gets 1.0e-8$, max epochs $\gets 300$}
\State $\bm{\tau}^0 = \bm{0}_n, \ \hat{\bm{\tau}}^0 = \bm{0}_p, \bm{\xi}^{0} = \mathbf{X}\hat{\bbeta}^{0}$
\State flag $=$ True,  err $= \infty$, t $= 0$ 
\While{err $\geq$ tol and flag}
\State $t \gets t + 1$
\State $\xi_{i}^{t} \gets \mathbf{X}_{i}'\hat{\bbeta}^{t-1} + \tau_i \dot{\mathcal{M}}_{g(. , T)}(\xi^{t-1}, \tau^{t-1}_i),\ \forall i$
\State $1/\hat{\tau}_{\mu}^{t} \gets \sum_{i=1}^n X^2_{i,\mu} \ddot{\mathcal{M}}_{g(. , T)}(\xi^t_i, \tau^{t-1}_i), \ \forall \mu $ %\frac{\ddot{g}(\xi_i^{t}, T)}{1 + \tau_i^{t} \ddot{g}(\xi_i^{t}, T)}
\State $\psi_{\mu}^t \gets  \hat{\beta}_{\mu}^{t-1} - \hat{\tau}_{\mu}^t\sum_{i=1}^n \dot{\mathcal{M}}_{g(. , T)}(\xi^t_i, \tau^{t-1}_i)X_{i,\mu}, \ \forall \mu$
\State $\hat{\beta}_{\mu}^{t}\gets \prox_{r(.)}(\psi_{\mu}^t, \hat{\tau}_{\mu}^t), \ \forall \mu $
\State $\tau_i^{t} \gets \sum_{\alpha = 1}^p   \hat{\tau}_{\alpha}^t \ \prox'_{r(.)}(\psi_{\alpha}^{t}, \hat{\tau}_{\alpha}^t) X_{i,\alpha}^2, \ \forall i$
\State err $\gets \sqrt{ \|\hat{\bbeta}^{t} - \hat{\bbeta}^{t-1}\|_2^2 + \|\bm{\tau}^t - \bm{\tau}^{t-1}\|_2^2 +  \|\bm{\xi}^t - \bm{\xi}^{t-1}\|_2^2 + \|\hat{\bm{\tau}}^t - \hat{\bm{\tau}}^{t-1}\|_2^2}$
\If{t $\geq$ max epochs}
\State{flag $=$ False}
\EndIf
\EndWhile
\end{algorithmic}
\end{algorithm}

This can be efficiently implemented. It is advised to damp the AMP iteration to improve numerical stability and we notice empirically that this is the case also here.

\subsection{Approximate Message Passing for i.i.d. Gaussian covariates}
When 
\begin{equation}
    \mathbf{X} \sim \mathcal{N}(\bm{0}_p, \bm{I}_p/p)  \  ,
\end{equation}
one can further simplify the algorithm and obtain the Approximate Message Passing algorithm originally proposed to solve the Lasso optimization problem \cite{Maleki_2010} and then extended to generalized linear models in \cite{Rangan_2010}.
The algorithm can be stated elegantly in vector notation, which provides also a more compact formulation. To this end it is useful to introduce the average operator$\langle\rangle$, which for a vector $\bx\in \mathbb{R}^l$ reads 
\begin{equation}
    \langle \mathbf{x}\rangle = \frac{1}{l} \sum_{k=1}^l x_k \ .
\end{equation}
\begin{algorithm}
\caption{AMP algorithm for MAP in GLMs (vector notation)}\label{alg:cap}
\begin{algorithmic}
\Require{$\hat{\bbeta}^{0} \gets \bm{0}_p$,  tol $\gets 1.0e-8$, max epochs $\gets 300$}
\State $\bm{\tau}^0 = 0, \ \hat{\bm{\tau}}^0 = 0, \bm{\xi}^{0} = \mathbf{X}\hat{\bbeta}^{0}$
\State flag $=$ True,  err $= \infty$, t $= 0$ 
\While{err $\geq$ tol and flag}
\State $t \gets t + 1$
\State $\bm{\xi}^{t} \gets \mathbf{X}\hat{\bbeta}^{t-1} + \tau \dot{\mathcal{M}}_{g(. , T)}(\bm{\xi}^{t-1}, \tau^{t-1})$
\State $\hat{\tau}^{t} \gets \zeta / \Big\langle \ddot{\mathcal{M}}_{g(. , T)}(\bm{\xi}^t, \tau^{t-1}) \Big\rangle $ 
\State $\bm{\psi}^t \gets  \hat{\bbeta}^{t-1} - \hat{\tau}^t \mathbf{X}'\dot{\mathcal{M}}_{g(. , T)}(\bm{\xi}^t, \tau^{t-1})$
\State $\hat{\bbeta}^{t}\gets \prox_{r(.)}(\bm{\psi}^t, \hat{\tau}^t)$
\State $\tau^{t} \gets    \hat{\tau}^t  \ \Big\langle \prox'_{r(.)}(\bm{\psi}^{t},\hat{\tau}^t) \Big\rangle$
\State err $\gets \sqrt{ \|\hat{\bbeta}^{t} - \hat{\bbeta}^{t-1}\|_2^2 + (\tau^t - \tau^{t-1})^2 +  \|\bm{\xi}^t - \bm{\xi}^{t-1}\|_2^2 + (\hat{\tau}^t - \hat{\tau}^{t-1})^2}$
\If{t $\geq$ max epochs}
\State{flag $=$ False}
\EndIf
\EndWhile
\end{algorithmic}
\end{algorithm}

\newpage
\section{COX-AMP : an AMP algorithm for the Cox model}
The optimization of the Cox Partial Likelihood can be regarded as an alternating minimization, as we have already argued in the derivation of the Coordinate-wise Descent algorithm, see \ref{app : cd_algorithm}. At the fixed point we must have that 
\begin{eqnarray}
    \dot{g}(\mathbf{X}\hat{\bbeta}_n , \hat{\Lambda}_n(\bT), \bm{\Delta}) = \bm{0}\\
    \hat{\Lambda}_n(T_i) = {\rm NA}(T_i, \mathbf{X}\hat{\bbeta}_n) \ .
\end{eqnarray}
The idea now is to do a step with Generalized - AMP in $\bbeta$ at fixed $\Lambda$ and then update $\Lambda$ with the Nelson Aalen estimator (\ref{app: NA_est}) as 
\begin{equation}
    \hat{\Lambda}_n^{t}(\bT) \gets {\rm NA}\big(\bT, \prox_{g(. \hat{\Lambda}_n^{t-1}(\bT), \bm{\Delta})}(\bm{\xi}^t, \tau^{t-1})\big)
\end{equation}
since at the fixed point of the Generalized -AMP algorithm, we must have that 
\begin{equation}
    \prox_{g(. \hat{\Lambda}_n^t(\bT), \bm{\Delta})}(\bm{\xi}, \tau) = \mathbf{X}\hat{\bbeta}_n  \ .
\end{equation}
The algorithm is schematized in the main text, see Algorithm \ref{alg : amp_cox}.

\end{document}